\documentclass[preprint,12pt]{elsarticle}

\newcommand{\sech}{\operatorname{sech}}
\newcommand{\sgn}{\operatorname{sgn}}

\usepackage{amsmath}
\usepackage{amssymb}
\usepackage{amsthm}
\usepackage{subfigure}
\journal{}%lsJournal of Computational Physics}

\begin{document}

\begin{frontmatter}

\author[mata]{Francisco de la Hoz}
\ead{francisco.delahoz@ehu.eus}
\address[mata]{Department of Applied Mathematics and Statistics and Operations Research, Faculty of Science and
Technology, University of the Basque Country UPV/EHU, Barrio Sarriena S/N, 48940 Leioa, Spain}

\author[matb]{Carlota M. Cuesta\corref{cor1}}
\ead{carlotamaria.cuesta@ehu.eus}
\address[matb]{Department of Mathematics, Faculty of Science and Technology, University of the
Basque Country UPV/EHU, Barrio Sarriena S/N, 48940 Leioa, Spain}

\cortext[cor1]{Corresponding author}

\title{A pseudo-spectral method for a non-local KdV-Burgers equation posed on $\mathbb R$}

\begin{abstract}

In this paper, we present a new pseudo-spectral method to solve the initial value problem associated to a non-local KdV-Burgers equation involving a Caputo-type fractional derivative. The basic idea is, using an algebraic map, to transform the whole real line into a bounded interval where we can apply a Fourier expansion. Special attention is given to the correct computation of the fractional derivative in this setting.

\end{abstract}


\end{frontmatter}


\section{Introduction}

\label{s:introduction}

In this paper, we present a new pseudo-spectral method to solve the initial value problem associated to the following non-local KdV-Burgers equation \cite{ACH}:
\begin{equation}
\label{e:uoverR}
\partial_t v + \partial_x(v^2) = \partial_x\mathcal D^\alpha v + \tau\partial_x^3 v, \quad x\in\mathbb R,\ t\ge 0,
\end{equation}

\noindent with $\tau > 0$, where $\mathcal D^\alpha$ denotes the non-local operator, acting only on the variable $x\in \mathbb{R}$,
\begin{equation}
\label{e:Da}
\mathcal D^\alpha v(x) = d_\alpha\int_{-\infty}^x\frac{v'(y)}{(x-y)^\alpha}dy,
\end{equation}

\noindent with $0 < \alpha < 1$, and
\begin{equation*}
d_\alpha = \frac{1}{\Gamma(1 - \alpha)} > 0,
\end{equation*}

\noindent being $\Gamma$ the Gamma function. The operator $\mathcal D^\alpha$ can be regarded as a left-sided fractional derivative in the Caputo sense, see e.g. \cite{Kilbas}, with integration taken from $-\infty$.

Equation (\ref{e:uoverR}) with $\alpha=1/3$ and either a quadratic flux, as above, or a cubic one, has been derived from a model of two-layer shallow water flow. This equation results by performing formal asymptotic expansions associated to the triple-deck (boundary
layer) theory used in fluid mechanics (see, e.g. \cite{Kluwick2010} and \cite{NV}). It is well-known that hyperbolic systems, such as the two dimensional shallow water model, may exhibit solutions that develop discontinuities or shocks. These shocks correspond to the so-called hydraulic-jumps or bores in this context. In order to investigate the internal structure of shock waves, it is customary to consider the viscosity of the system or to introduce it artificially. In the particular case studied in \cite{Kluwick2010}, viscosity terms become important near the bottom boundary. In order to take its effect into account, a boundary layer needs to be introduced. Additionally, the effect due to the stream curvature that becomes important near shock waves cannot be neglected.

In \cite{Kluwick2010} (for the single layer case) and in \cite{NV} (see also \cite{KSC}), the authors keep track of the viscosity on the bottom boundary in a way consistent with the Navier-Stokes equations, by means of matched asymptotic expansions in the limit of the Reynolds number tending to infinity. As a consequence, a distinguished limit associated with a \emph{triple deck} problem is presented. The main idea of the triple deck theory (see e.g. \cite{29}) is that the leading order outer flow is not independent of the boundary layer flow, unlike in the classical Prandtl's theory \cite{28}. These considerations, in a weak interaction limit, lead to an equation of the type (\ref{e:uoverR}), where the fractional derivative term results from the effect of viscosity, and the third order term results from the stream curvature. The infinite domain where the equation is posed corresponds to the domain of the inner region in the asymptotic limit. In this region, traveling wave solutions and their stability are naturally analyzed, since they resemble the inner structure of shock waves.

There are other models where such a fractional derivative terms arise (see v.g. \cite{chester}, \cite{keller}, \cite{sugimoto89}, \cite{sugimoto90} and \cite{sugiJFM}). The common feature of these models is that they result from the analysis of a boundary layer and thus they need to be studied in an unbounded domain.

The existence of the Cauchy problem for (\ref{e:uoverR}), and the existence of traveling wave solutions and their stability have been studied rigorously in \cite{AHS} (with $\tau=0$), and in \cite{ACH} (with $\tau>0$). It is worth emphasizing that the study of traveling wave solutions requires naturally that the equation is posed on $\mathbb R$, since these waves travel along the domain and are expected to emerge in a large time limit. Our aim is thus to present a numerical method that deals accurately and efficiently with the problem posed on $\mathbb R$, taking particular care of the numerical treatment of ${\cal D}^{\alpha}v$.

In the following pages, we will develop a new pseudo-spectral method for \eqref{e:uoverR}. Essential references on spectral methods can be found in \cite{ShenTangWang2011,stenger2011,Canutoetal2007,boyd2001,fornberg,trefethen}, together with the more classical \cite{Canuto1988,gottlieb}.

One of the main difficulties in dealing with \eqref{e:uoverR} is the unboundedness of the spatial domain. However, according to Boyd \cite[p. 338]{boyd2001}, the many possible options for unbounded domains always fall into one of three following broader categories:
\begin{enumerate}

\item Truncating the domain (taking $x\in[-L, L]$, with $L\gg 1$).

\item Using basis functions intrinsic to an infinite interval (for example, Hermite functions, sinc functions).

\item Mapping the unbounded interval to a finite interval, followed by the application of Chebyshev polynomials or of a Fourier series.

\end{enumerate}

We will adopt the third option here. It is possible to generate a great variety of new basis sets for the infinite interval, which are the images of Chebyshev polynomials or Fourier series \cite[p. 355]{boyd2001} under a change of the independent variable $x$ that maps $\mathbb{R}$ into a finite interval. Although an infinite variety of maps is possible, we will concentrate on a very important one, the so-called algebraic map (see v.g. \cite{GroschOrszag1977} and \cite{Boyd1987}):
\begin{equation}
\label{e:algebraicmap}
\xi = \dfrac{x}{\sqrt{L^2 + x^2}} \Longleftrightarrow x = \dfrac{L\xi}{\sqrt{1 - \xi^2}},
\end{equation}

\noindent with $L > 0$, which maps the whole real line $x\in\mathbb R$ into the interval $\xi\in[-1, 1]$ and vice versa. Then, we can use the Chebyshev polynomials (of the first kind) $T_k(\xi)$ over the new domain $\xi\in[-1,1]$:
\begin{equation}
\label{e:Tk}
T_k(\xi) = \cos(k\arccos(\xi)), \quad \forall k\in\mathbb N\cup\{0\}.
\end{equation}

\noindent Moreover, \eqref{e:algebraicmap} allows us to define a basis set for the infinite interval, formed by the so-called rational Chebyshev polynomials:
\begin{equation}
\label{e:TBk}
TB_k(x) = T_k\left(\dfrac{x}{\sqrt{L^2+x^2}}\right), \quad x\in\mathbb R, \quad \forall k\in\mathbb N\cup\{0\}.
\end{equation}

\noindent These functions form an orthogonal basis in $\mathbb R$ with respect to the weight $(1+ x^2)^{-1}$:
\begin{equation*}
\int_{-\infty}^{+\infty}\dfrac{TB_m(x)TB_n(x)}{1 + x^2}dx =
\begin{cases}
\pi / 2, & m = n > 0,
\\
\pi, & m = n = 0,
\\
0, & m \not= n.
\end{cases}
\end{equation*}

\noindent Rational Chebyshev polynomials appear to behave extremely well in a great variety of problems. For instance, in \cite{delahozvadillo}, rational Chebyshev polynomials, Hermite functions and sinc functions were compared when solving numerically two and three-dimensional nonlinear diffusion equations over unbounded domains. While Hermite and sinc functions are adequate for the approximation of functions with exponential decay, it was shown in \cite{delahozvadillo} that rational Chebyshev polynomials are the most versatile option, because, besides being a good choice for approximating exponentially decaying functions, they really excel when applied to polynomially decaying solutions. In fact, in this paper, we consider an exponentially decaying initial data for \eqref{e:uoverR} that become algebraically decaying as $t$ increases (see \cite{ACH}). Hence, in our opinion, rational Chebyshev polynomials are here the best option, and, thus, we have preferred them over other choices like Hermite functions and sinc functions.

In general, an advantage of not truncating the domain is that the boundary conditions can often be ignored when the domain of integration is infinite (see \cite{Boyd1987}), while truncating the domain necessitates setting artificial boundary conditions (see v.g. \cite{tsynkov1998} for a very complete review on artificial boundary conditions in domain truncation problems).

Although we can work directly with \eqref{e:TBk}, it is easier to use them under a trigonometric representation, i.e., through the change of variable
\begin{equation}
\label{e:changeofvariable}
\xi = \cos(s) \Longleftrightarrow x = L\cot(s), \quad s\in[0, \pi].
\end{equation}

\noindent Then, \eqref{e:Tk} and \eqref{e:TBk} become
\begin{equation*}
T_k(\xi) = TB_k(x) = \cos(ks),
\end{equation*}

\noindent i.e., we reduce the development of a function $v(x)$ over $\mathbb R$ in a series of rational Chebyshev polynomials to obtaining just a cosine expansion:
\begin{equation*}
u(s) \equiv v(x(s)) = \sum_{k = 0}^\infty \hat a(k)\cos(ks).
\end{equation*}

\noindent We recall that, after applying \eqref{e:changeofvariable}, $s\in[0, \pi]$, while the cosines are periodic in $s\in[0, 2\pi]$. Therefore, in order to compute the cosine expansion of a function defined over half a period, we may use an even extension of that function. At this point, it is very important to underline that an even extension involving only cosines is by no means the only option. Indeed, we will consider a more general series expansion of $v(x)$ under the change of variable $x = L\cot(s)$:
\begin{equation}
\label{e:ucossin}
u(s) = \sum_{k = 0}^\infty \hat a(k)\cos(ks) + \sum_{k = 1}^\infty \hat b(k)\sin(ks),
\end{equation}

\noindent which implies augmenting \eqref{e:TBk} by a second set of basis functions denoted $SB_k(x)$ and defined \cite[p. 127]{Boyd1987} by:
\begin{equation*}
SB_{k+1}(x) = (1-\xi^2)^{1/2}U_k(\xi) = \sin((k+1)s), \quad \forall k\in\mathbb N\cup\{0\},
\end{equation*}

\noindent where $U_k(\xi)$ are the so-called Chebyshev polynomials of the second kind:
\begin{equation*}
U_k(\cos(s)) = \frac{\sin((k+1)s)}{\sin(s)}, \quad \forall k\in\mathbb N\cup\{0\}.
\end{equation*}

\noindent Even if we could work directly with $TB_k(x)$ and $TS_k(x)$, we follow the strategy recommended by Boyd \cite{Boyd1987}, i.e., to change to the trigonometric variable and then apply Fourier series.

The structure of this paper is as follows. In Section \ref{s:method}, we propose a pseudo-spectral method for simulating \eqref{e:uoverR}, which is based on the change of variable \eqref{e:changeofvariable}, together with the application of a Fourier series expansion. Although there are many references on the numerical computation by spectral methods of fractional derivatives, most of them do not seem to consider integration from $-\infty$ or to $+\infty$, as is necessary in \eqref{e:Da}. For instance, \cite{miyakoda2009} (which is based on \cite{hasegawa1991}) and \cite{bhrawy2013} (for the closely related fractional integration) make use of the Chebyshev polynomials, and in \cite{Li2012}, Legrende, Chebyshev and Jacobi polynomials are used for the computation of the Caputo derivative. Other spectral methods use a suitable generalization of Jacobi polynomials for the computation of fractional derivatives (see for instance \cite{zayernouri2014} and the references therein). We remark that these generalized Jacobi polynomials are well suited in the case of bounded domains; since, in that case, they form an orthogonal basis, and the elements of this basis are mapped into the basis by application of fractional derivatives. Fractional derivatives are also implemented in the \textsc{Chebfun} package from \textsc{Matlab\copyright}, where the Gauss-Jacobi quadrature is used (see v.g. \cite{hale2012}). In \cite{Huang2014}, a spectral method for substantial fractional differential equations on a semi-infinite domain is developed, without truncating the domain; notice that this fractional derivative has a much more regular kernel. In any case, we have found no references on the computation of the fractional derivative of rational Chebyshev polynomials. Hence, the central part of this paper is the accurate computation of the operator \eqref{e:Da}, or, more precisely, of $\partial_x\mathcal D^\alpha$, which is done in Section \ref{s:fractional}; whereas, in Section \ref{s:regularity}, we do a numerical study on the minimum regularity required by the involved functions.

The structure of \eqref{e:uoverR} strongly suggests the use of an implicit-explicit (IMEX) scheme in time, where the highest-order term, $\tau\partial_x^3$, appears implicitly, and the other lower-order terms appear explicitly; this is done in Section \ref{s:timediscretization}. The numerical experiments are carried out in Section \ref{s:experiments}, and show that we can simulate \eqref{e:uoverR} until extremely large times.

Let us finish this introduction by mentioning that, in recent years, there has been an increasing interest in solving Fractional Partial Differential Equations; in particular, equations with non-local diffusions, as in \eqref{e:uoverR}, (see for instance \cite{buenoorovio2014} for fractional Laplacian); \cite{tian2014} for advection-diffusion with non-local diffusion; and \cite{Zeng2014}, for a spectral method applied to non-local reaction-diffusion equations). Moreover, a computation of the fractional evolution Burgers equation can be found in \cite{Zhao2015}.

\section{Numerical method}

\label{s:method}

As explained above, in order to simulate numerically \eqref{e:uoverR}, we map the space domain from $\mathbb R$ into $[0,\pi]$ via the change of variable $x = L\cot(s)$, and then apply Fourier series to the resulting equation. We discretize the space variable at the non-final nodes:
\begin{equation}
\label{e:nodes}
s_j = \frac{\pi(2j+1)}{2N} \Longleftrightarrow  x_j = L\cot\left(\frac{\pi(2j+1)}{2N}\right) \Longleftrightarrow  \xi_j = \cos\left(\frac{\pi(2j+1)}{2N}\right),
\end{equation}

\noindent where $0 \le j \le N-1$, which divide the interval $[0, \pi]$ in $N$ equally-spaced parts. Again, $u(s_j)\equiv v(x(s_j)) = v(x_j)$.

As mentioned in the introduction, instead of working only with cosines, we consider a more general series expansion of $u$ in the form \eqref{e:ucossin}. However, the implementation is done with more ease by expanding $u$ into series of $e^{iks}$:
\begin{equation}
\label{e:uexp}
u(s) = \sum_{k = -\infty}^\infty \hat u(k)e^{iks}, \quad s\in[0, \pi],
\end{equation}

\noindent where, in order to determine the coefficients $\hat u(k)$, it is necessary to decide how to extend the function from $s\in[0, \pi]$ to $s\in[0, 2\pi]$; i.e., we discretize $[0, 2\pi]$ in $2N$ nodes, $s_j = \pi(2j + 1) / (2N)$, $0\le j \le 2N-1$. We recall that, if we consider an even extension, we have a cosine expansion, while an odd extension yields a sine expansion.

In general, given a $2N$-term approximation of a $2\pi$-periodic function $u$
\begin{equation}
\label{e:approxu(s)}
u(s) \approx \sum_{k = -N}^{N-1} \hat u(k)e^{iks}, \quad s\in[0, 2\pi],
\end{equation}

\noindent it is possible to determine uniquely the $\hat u(k)$ by evaluating \eqref{e:approxu(s)} at $s_j$
\begin{equation}
\label{e:idft}
u(s_j) = \sum_{k = -N}^{N-1} \hat u(k)e^{ik\pi(2j + 1) / (2N)}  = \sum_{k = 0}^{2N-1} [\hat u(k)e^{ik\pi / (2N)}]   e^{2ijk\pi / (2N)};
\end{equation}

\noindent then,
\begin{equation}
\label{e:dft}
\hat u(k) = \frac{e^{-ik\pi / (2N)}}{2N}\sum_{j = 0}^{2N-1} u(s_j)e^{-2ijk\pi / (2N)},
\end{equation}

\noindent where $0\le i, j\le 2N-1$. Furthermore, the discrete Fourier transforms \eqref{e:idft} and \eqref{e:dft} can be computed very efficiently by means of the fast Fourier transform (FFT) \cite{FFT}. Observe that, although \eqref{e:approxu(s)} is in principle exact only at the $s_j$, it is nonetheless a very good approximation for all $s$, provided that $u$ is regular enough and that $N$ is large enough. Let us mention also that, in order to \emph{clean the spectrum}, it is convenient to round to zero those $\hat u(k)$ whose absolute value is smaller than the epsilon of the machine, $\varepsilon = 2.2204\cdot10^{-16}$.

Let us express now all the terms of \eqref{e:uoverR} in function of the new space-variable $s$. It is straightforward \cite{Boyd1987} to check that
\begin{equation}
\label{e:u_xu_s}
\begin{split}
 v_x & = -\frac{\sin^2(s)}{L}u_s,
    \cr
 v_{xx} & = \frac{\sin^4(s)}{L^2}u_{ss} + \frac{2\sin^3(s)\cos(s)}{L^2}u_s,
    \cr
 v_{xxx} & = - \frac{\sin^6(s)}{L^3}u_{sss} - \frac{6\sin^5(s)\cos(s)}{L^3}u_{ss}
    + \frac{\sin^4(s)(8\sin^2(s) - 6)}{L^3}u_s.
\end{split}
\end{equation}

\noindent Therefore, given $u$ in the form of \eqref{e:approxu(s)}, we can compute $v_x$, $v_{xx}$ and $v_{xxx}$ very efficiently. Indeed, let us suppose that $u$ consists of just one node, i.e., $u = e^{iks}$. Then, representing $\cos(s)$ and $\sin(s)$ in polar form, \eqref{e:u_xu_s} becomes
\begin{equation}
\label{e:udiffexp}
\begin{split}
 v_x = & \, \frac{ik}{4L}e^{i(k+2)s} - \frac{ik}{2L}e^{iks} + \frac{ik}{4L}e^{i(k-2)s},
    \cr
 v_{xx} = & - \frac{k^2+2k}{16L^2}e^{i(k+4)s} + \frac{k^2+k}{4L^2}e^{i(k+2)s} - \frac{3k^2}{8L^2}e^{iks}
    \cr
& + \frac{k^2-k}{4L^2}e^{i(k-2)s} - \frac{k^2-2k}{16L^2}e^{i(k-4)s},
    \cr
 v_{xxx} = & -i\frac{k^3+6k^2+8k}{64L^3}e^{i(k+6)s} + i\frac{3k^3+12k^2+12k}{32L^3}e^{i(k+4)s}
    \cr
& - i\frac{15k^3+30k^2+24k}{64L^3}e^{i(k+2)s} + i\frac{5k^3+4k}{16L^3}e^{iks}
    \cr
& - i\frac{15k^3-30k^2+24k}{64L^3}e^{i(k-2)s} + i\frac{3k^3-12k^2+12k}{32L^3}e^{i(k-4)s}
    \cr
& - i\frac{k^3-6k^2+8k}{64L^3}e^{i(k-6)s}.
\end{split}
\end{equation}

\noindent From \eqref{e:udiffexp}, we can construct immediately (very sparse) differentiation matrices in the Fourier side.

At this point, some comment about the choice of $L$ is required; indeed, the correct choice of the scaling is a matter of concern. Even if there are some theoretical results \cite{Boyd1982}, the optimal scaling values depend on the number of points, the class of functions, and the type of problem, and can even change during time. However, a good working rule of thumb seems to be that the absolute value of the function at the extremal grid points is smaller than an accuracy threshold $\varepsilon$.

\subsection{Computation of the fractional derivative}

\label{s:fractional}

The most involved part in this paper is the accurate computation of the fractional derivative $\partial_x\mathcal D^\alpha v(x)$. The operator $\mathcal D^\alpha v(x)$, in the form of \eqref{e:Da}, is not very well-suited for numerical calculations when $\alpha$ is close to one, because $1 / \Gamma(1 - \alpha)$ tends to infinity as $\alpha\to1^{-}$. On the other hand, throughout this paper, we are considering solutions that remain at least $\mathcal C^3$ for all time (see \cite{ACH}). Therefore, we can safely integrate \eqref{e:Da} by parts:
\begin{align}
\label{e:byparts}
\mathcal D^\alpha v(x) & = \frac{1}{\Gamma(1-\alpha)}\lim_{\substack{b\to x^- \\ a\to-\infty}}\int_a^b
v'(y)(x-y)^{-\alpha}dy
    \cr
& = \frac{-1}{\Gamma(2-\alpha)}\lim_{\substack{b\to x^- \\ a\to-\infty}}\left[v'(y)(x-y)^{1-\alpha}\Big|_{y=a}^{y=b} - \int_a^b
v''(y)(x-y)^{1-\alpha}dy\right]
    \cr
& = \frac{1}{\Gamma(2-\alpha)}\int_{-\infty}^x v''(y)(x-y)^{1-\alpha}dy,
\end{align}

\noindent because, for a fixed $x$,
\begin{equation}
\label{e:decay1}
\lim_{y\to x^-}v'(y)(x - y)^{1-\alpha} = 0, \qquad
\lim_{y\to -\infty}v'(y)(x - y)^{1-\alpha} = 0.
\end{equation}

\noindent The last condition is satisfied by all the functions we are interested in; for instance, in \eqref{e:uoverR}, one expects that $v(y,t)$ tends exponentially to a constant as $y\to-\infty$, if the initial condition decays sufficiently fast (see \cite{ACH}). A similar reasoning enables us to conclude that $\partial_x\mathcal D^\alpha v(x)$ and $\mathcal D^\alpha v'(x)$ are equivalent. Indeed, differentiating \eqref{e:byparts}:
\begin{align}
\label{e:DxDa}
\partial_x\mathcal D^\alpha v(x) & = \frac{1}{\Gamma(2-\alpha)}\partial_x\left(\int_{-\infty}^x v''(y)(x-y)^{1-\alpha}dy\right)
    \cr
& = \frac{1}{\Gamma(2-\alpha)}\bigg(v''(y)(x-y)^{1-\alpha}\Big|_{y=x} + \int_{-\infty}^x \partial_x\left(v''(y)(x-y)^{1-\alpha}\right)  dy\bigg)
    \cr
& = \frac{1}{\Gamma(1-\alpha)}\int_{-\infty}^x v''(y)(x-y)^{-\alpha} dy
    \cr
& = \frac{1}{\Gamma(2-\alpha)}\int_{-\infty}^x v'''(y)(x-y)^{1-\alpha} dy,
\end{align}

\noindent where we are assuming that, for a fixed $x$,
\begin{equation}
\label{e:decay2}
\lim_{y\to -\infty}v''(y)(x - y)^{1-\alpha} = 0.
\end{equation}

\noindent In the numerical method that we propose, it is compulsory to integrate \eqref{e:Da} once by parts, as is done in \eqref{e:byparts}. On the other hand, instead of working with \eqref{e:DxDa}, it is possible to first compute numerically \eqref{e:byparts}, and then differentiate it numerically. However, we have chosen to work with \eqref{e:DxDa}, because it is more compact to implement and more accurate.

Let us consider the change of variable $x = L\cot(s)$ from \eqref{e:changeofvariable}; then, using \eqref{e:u_xu_s},
\begin{align}
v_{yyy}dy = \left[\frac{\sin^4(\eta)}{L^2}u_{\eta \eta \eta} + \frac{6\sin^3(\eta)\cos(\eta)}{L^2}u_{\eta \eta}
    - \frac{\sin^2(\eta)(8\sin^2(\eta) - 6)}{L^2}u_\eta\right]d\eta.
\end{align}

\noindent Hence, \eqref{e:DxDa} becomes
\begin{align}
\label{e:intspi}
\left(\partial_x\mathcal D^\alpha v\right)(x(s)) = & -L^{-1-\alpha}\int_s^{\pi}\big[\sin^4(\eta)u'''(\eta) + 6\sin^3(\eta)\cos(\eta)u''(\eta)
    \cr
& - \sin^2(\eta)(8\sin^2(\eta) - 6)u'(\eta)\big](\cot(s) - \cot(\eta))^{1-\alpha}d\eta
    \cr
= & \int_s^{\pi}w(\eta)(\cot(s) - \cot(\eta))^{1-\alpha}d\eta,
\end{align}

\noindent where
\begin{align}
\label{e:w}
w(s) \equiv\, &  L^{-1-\alpha}[-\sin^4(s)u'''(s) - 6\sin^3(s)\cos(s)u''(s)
    \cr
& + \sin^2(s)(8\sin^2(s) - 6)u'(s)].
\end{align}

\noindent We want to approximate \eqref{e:intspi} at the nodes $s_j$ defined in \eqref{e:nodes}. At this point, instead of applying a quadrature formula like the trapezoidal rule directly to \eqref{e:intspi}, we have chosen to transform it into an integral over $[0, \pi]$ by means of a characteristic function:
\begin{equation}
\label{e:int0pi}
\partial_x\mathcal D^\alpha u(s) \equiv \left(\partial_x\mathcal D^\alpha v\right)(x(s)) = \int_0^{\pi}w(\eta)\chi_{[s,\pi]}(\eta)(\cot(s) - \cot(\eta))^{1-\alpha}d\eta.
\end{equation}

\noindent Now, given a function $f$ defined over $[0, \pi]$ and whose value is known over the nodes \eqref{e:nodes}, a way of approximating its integral is by means of the well-known Chebyshev-Gauss quadrature \cite{abramowitz}:
\begin{equation}
\label{e:quadrature}
\int_0^\pi f(\eta)d\eta = \int_{-1}^{1} \frac{f(\arccos(\xi))}{\sqrt{1-\xi^2}}d\xi \approx \frac{\pi}{N}\sum_{j=0}^{N-1}f(\arccos(\xi_j)) = \frac{\pi}{N}\sum_{j=0}^{N-1}f(s_j).
\end{equation}

\noindent However, in practice, we have found that applying this last formula to \eqref{e:int0pi} gives poor results; because, inside the integrand, $\chi_{[s_j,\pi]}(\eta)(\cot(s_j) - \cot(\eta))^{1-\alpha}$ has the singularity point precisely at $\eta = s_j$. Therefore, we consider instead the families of nodes
\begin{equation*}
s_j^{(m)} = \frac{\pi(2j+1)}{2^{m+1}N}, \quad 0 \le j \le 2^{m}N-1, \quad m = 1, 2, \ldots
\end{equation*}

\noindent Since $s_j^{(m)}$ never coincides with any $s_j$, we avoid evaluating $\chi_{[s_j,\pi]}(\eta)(\cot(s_j) - \cot(\eta))^{1-\alpha}$ at the singularity point, so we can apply \eqref{e:quadrature} to \eqref{e:int0pi} at the new nodes, to obtain an approximation that depends on $m$:
\begin{align}
\label{e:int0pim}
\left[\partial_x\mathcal D^\alpha\right]^{(m)} u(s_j) & \equiv \left[\partial_x\mathcal D^\alpha\right]^{(m)} v(x_j)
    \cr
& \equiv \frac{\pi}{2^mN}\sum_{l = 0}^{2^mN-1} w(s_l^{(m)})\chi_{[s_j,\pi]}(s_l^{(m)})(\cot(s_j) - \cot(s_l^{(m)}))^{1-\alpha}
    \cr
& = \frac{\pi}{2^mN}\sum_{l = 2^{m-1}(2j+1)}^{2^mN-1} w(s_l^{(m)})(\cot(s_j) - \cot(s_l^{(m)}))^{1-\alpha},
\end{align}

\noindent where, in the last line, we have used that
\begin{align}
s_l^{(m)} \ge s_j \Longleftrightarrow \frac{\pi(2l+1)}{2^{m+1}N} \ge \frac{\pi(2j+1)}{2N} \Longleftrightarrow
l \ge 2^{m-1}(2j+1).
\end{align}

\noindent Bearing in mind the Fourier decomposition \eqref{e:approxu(s)} of $u$, it is possible to evaluate \eqref{e:w} at $s_l^{(m)}$ fast and accurately. Indeed, if $u = e^{ik s}$, applying it to \eqref{e:w}, we get an expression similar to those in \eqref{e:udiffexp}:
\begin{equation}
\label{e:wexp}
\begin{split}
w_k(s) =\ & i\frac{k^3+6k^2+8k}{16L^{1+\alpha}}e^{i(k+4)s} - i\frac{k^3+3k^2+2k}{4L^{1+\alpha}}e^{i(k+2)s} + i\frac{3k^3}{8L^{1+\alpha}}e^{iks}
    \cr
& - i\frac{k^3-3k^2+2k}{4L^{1+\alpha}}e^{i(k-2)s} + i\frac{k^3-6k^2+8k}{16L^{1+\alpha}}e^{i(k-4)s};
\end{split}
\end{equation}

\noindent where we need to consider only $k > 0$, because $w_{-k}(s) = \bar w_{k}(s)$. Moreover, from a computational point of view, it is much faster to evaluate $w_k(s)$ if rewritten as
\begin{equation}
\label{e:wexp2}
\begin{split}
w_k(s) = \frac{e^{iks}}{L^{1+\alpha}}\bigg[& i\frac{k^3+8k}{8}\cos(4s) - \frac{3k^2}{4}\sin(4s)
    \cr
& - i\frac{k^3+2k}{2}\cos(2s) + \frac{3k^2}{2}\sin(2s) + i\frac{3k^3}{8}\bigg].
\end{split}
\end{equation}

\noindent In \eqref{e:int0pim}, we are considering only $s\in[0, \pi]$, although we are working with $s\in[0,2\pi]$. However, since $\cot(s)$ is $\pi$-periodic, we simply define
\begin{align}
\label{e:intpi2pim}
\left[\partial_x\mathcal D^\alpha\right]^{(m)} u(s_{j + N}) & \equiv \frac{\pi}{2^mN}\sum_{l = 2^{m-1}(2j+1)}^{2^mN-1} w(s_l^{(m)} + \pi)(\cot(s_j) - \cot(s_l^{(m)}))^{1-\alpha},
\end{align}

\noindent where $s_{j + N} = s_j + \pi$; moreover, if $u = e^{iks}$, it follows from \eqref{e:wexp} that
\begin{align}
\label{e:wpi}
w_k(s + \pi) = (-1)^kw_k(s).
\end{align}

\noindent Summarizing, given an arbitrary function $u(s)$ defined over $[0, \pi]$, we extend it to $[0, 2\pi]$; then we calculate its Fourier expansion \eqref{e:approxu(s)}; then we compute \eqref{e:int0pim} and \eqref{e:intpi2pim}, bearing in mind \eqref{e:wexp2} and \eqref{e:wpi}. For not extremely large $N$, the most efficient option, and the one we follow, is to generate an operational $(2N)\times(2N)$-matrix $\mathbf M_{\alpha, N}^{(m)}$ whose columns are precisely $\left[\partial_x\mathcal D^\alpha\right]^{(m)}e^{iks}$. Therefore, \eqref{e:int0pim} is reduced to just multiplying $\mathbf M_{\alpha, N}^{(m)}$ by the column vector formed by the Fourier coefficients \eqref{e:approxu(s)} of $u$, which we denote by $\widehat{\mathbf{U}}$:
\begin{equation}
\label{e:Dam}
\left[\partial_x\mathcal D^\alpha\right]^{(m)} u(s_j) \equiv \left[\mathbf M_{\alpha, N}^{(m)}\cdot\widehat{\mathbf{U}}\right]_j.
\end{equation}

\noindent After some optimization of the code, this matrix can be generated in a very efficient way; moreover, it needs to be calculated only once for given $\alpha$, $N$ and $m$, and then it can be stored. $L$ can be chosen equal to one and, for other choices of $L$, apply a scaling, i.e., take $L^{-1-\alpha}\mathbf M_{\alpha, N}^{(m)}$. Let us mention also that, if $u$ is real, it is convenient to impose this fact explicitly in \eqref{e:Dam}.

In order to test the quadrature formula \eqref{e:int0pim}, we have considered three functions satisfying \eqref{e:decay1} and \eqref{e:decay2}, and such that their fractional derivative $\partial_x\mathcal D^a$ can be explicitly computed by means of, for instance \textsc{Mathematica\copyright}. More precisely, we have considered a function with quadratic decay,
\begin{equation*}
v_1(x) = \frac{1}{1 + x^2},
\end{equation*}

\noindent such that
\begin{equation}
\label{e:dav1}
\begin{split}
\partial_x & \mathcal D^\alpha v_1(x) = -\frac{\pi\alpha (1 + \alpha) \csc(\alpha \pi)}{\Gamma(1 - \alpha)}(1 + x^2)^{-(3 + \alpha) / 2}
    \cr
& \cdot\Big[\sin\Big(\frac{\alpha \pi}{2} + (1 + \alpha)\arctan(x)\Big) + x \cos\Big(\frac{\alpha \pi}{2} + (1+\alpha)\arctan(x)\Big)\Big];
\end{split}
\end{equation}

\noindent a function with quartic decay,
\begin{equation*}
v_2(x) = \frac{1}{(1 + x^2)^2},
\end{equation*}

\noindent such that
\begin{equation}
\label{e:dav2}
\begin{split}
\partial_x & \mathcal D^\alpha  v_2(x) = \frac{\pi\alpha (1 + \alpha)}{4 \Gamma(1 - \alpha)} (1 + x^2)^{-3 - \alpha/2}
    \cr
& \cdot\Big[\sec\Big(\frac{\alpha\pi}{2}\Big)\big[ ((3\alpha + 8)x - \alpha x^3)\sin(\alpha \arctan(x))
    \cr
& + (-3 - \alpha + (6 + 3\alpha)x^2 + x^4)\cos(\alpha \arctan(x))\big]
    \cr
& + \csc\Big(\frac{\alpha\pi}{2}\Big)(1 + x^2)^{1/2}\big[
(-3-\alpha  + (1 + \alpha)x^2)   \sin((1 + \alpha) \arctan(x))
  \cr
& - (5x + 2 \alpha x + x^3) \cos((1 + \alpha) \arctan(x))\big]\Big];
\end{split}
\end{equation}

\noindent and, finally, a function with Gaussian decay,
\begin{equation*}
v_3(x) = \exp(-x^2),
\end{equation*}

\noindent such that
\begin{equation}
\label{e:dav3}
\begin{split}
\partial_x\mathcal D^\alpha v_3(x) = & \, \frac{1}{3\Gamma(1 - \alpha)}\bigg[ -
  6 (1 + \alpha) \Gamma\bigg(
    1 - \frac{\alpha}{2}\bigg) x{}_1F_1\bigg(\frac{3+\alpha}{2}, \frac{3}{2}, -x^2\bigg)
    \cr
& - 3\alpha \Gamma\bigg(\frac{1 - \alpha}{2}\bigg) {}_1F_1\bigg(1 + \frac{\alpha}{2}, \frac{3}{2}, -x^2\bigg)
    \cr
& + (4\alpha + 2\alpha^2) \Gamma\bigg(\frac{1 - \alpha}{2}\bigg) x^2 {}_1F_1\bigg(2 + \frac{\alpha}{2}, \frac{5}{2}, -x^2\bigg)\bigg],
\end{split}
\end{equation}

\noindent where ${}_1F_1(a, b, c)$ is the confluent hypergeometric function \cite{luke}, which can be accurately evaluated by means of, for instance, \textsc{Mathematica\copyright} or \textsc{Matlab\copyright}.

In general, it seems that the number of functions for which one can compute explicitly $\partial_x\mathcal D^\alpha$ is rather narrow. On the other hand, the resulting expressions are rather involved, as we can see in \eqref{e:dav1}, \eqref{e:dav2} and \eqref{e:dav3}, which thus constitute a stringent test for our numerical method. To measure the errors, we have evaluated \eqref{e:dav1}, \eqref{e:dav2}, and \eqref{e:dav3} for $0 < \alpha < 1$, which are the cases we are interested in; whereas, for the limiting cases $\alpha = 0$ and $\alpha = 1$, which can also be considered with our method, we have taken, respectively, the first and second derivatives of the test functions: $v_1'(x) = -2x(1 + x^2)^{-2}$, $v_1''(x) = (6x^2 - 2)(1 + x^2)^{-3}$; $v_2'(x) = -4x(1 + x^2)^{-3}$, $v_2''(x) = (20x^2 - 4)(1 + x^2)^{-4}$; and $v_3'(x) = -2x\exp(-x^2)$, $v_3''(x) = (4x^2-2)\exp(-x^2)$.

In our numerical experiments, we have taken $1001$ equally-spaced values of $\alpha\in[0,1]$, i.e., $\alpha = j / 1000$, $j = 0, \ldots, 1000$. After applying \eqref{e:changeofvariable}, and doing even extensions of the functions at $s = \pi$, we have multiplied the corresponding $\mathbf M_{\alpha, N}^{(m)}$ by their Fourier coefficients. The experiments corresponding to $v_1$ have been done with $L = 1.6$; those corresponding to $v_2$, with $L = 1.1$, and those corresponding to $v_3$, with $L = 4$. In all cases, $N = 64$.

In order to measure the accuracy of the results, we define, for a function $v(x)$ and a given $\alpha$, the following error related to the discrete $L^\infty$-norm:
\begin{equation}
\label{e:errorsEm}
E^{(m)}(\alpha) \equiv \max_j\left|\left[\mathbf M_{\alpha, N}^{(m)}\cdot\widehat{\mathbf U}\right]_j - \partial_x\mathcal D^\alpha v(x)\right|.
\end{equation}

\noindent For the sake of simplicity, we let $E_1^{(m)}(\alpha)$, $E_2^{(m)}(\alpha)$ and $E_3^{(m)}(\alpha)$ denote the errors corresponding to $v_1(x)$, $v_2(x)$ and $v_3(x)$, respectively. Then, Table \ref{t:testDam} shows $\max_{\alpha\in[0,1]}E_1^{(m)}(\alpha)$, $\max_{\alpha\in[0,1]}E_2^{(m)}(\alpha)$ and $\max_{\alpha\in[0,1]}E_3^{(m)}(\alpha)$, for $m = 1, \ldots, 6$. Even if the errors clearly decay as $m$ increases, they do it rather slowly.
\begin{table}[htb!]
\centering
\begin{tabular}{|c|c|c|c|}
\cline{2-4} \multicolumn{1}{c|}{} & $v_1(x)$ & $v_2(x)$ & $v_3(x)$
    \cr
\hline $\mathbf M^{(1)}$ &$5.0137\cdot10^{-3}$ &$8.4605\cdot10^{-3}$ &$1.2527\cdot10^{-2}$
    \cr
\hline $\mathbf M^{(2)}$ &$2.1906\cdot10^{-3}$ &$3.7336\cdot10^{-3}$ &$5.3076\cdot10^{-3}$
    \cr
\hline $\mathbf M^{(3)}$ &$9.7339\cdot10^{-4}$ &$1.6713\cdot10^{-3}$ &$2.3065\cdot10^{-3}$
    \cr
\hline $\mathbf M^{(4)}$ &$4.3810\cdot10^{-4}$ &$7.5653\cdot10^{-4}$ &$1.0206\cdot10^{-3}$
    \cr
\hline $\mathbf M^{(5)}$ &$1.9920\cdot10^{-4}$ &$3.4555\cdot10^{-4}$ &$4.5785\cdot10^{-4}$
    \cr
\hline $\mathbf M^{(6)}$ &$9.1329\cdot10^{-5}$ &$1.5902\cdot10^{-4}$ &$2.0763\cdot10^{-4}$
    \cr
\hline\end{tabular}
\caption{Given $v_1(x) = (1 + x^2)^{-1}$, $v_2(x) = (1 + x^2)^{-2}$, and $v_3 = \exp(-x^2)$, we have computed their fractional derivatives via $\mathbf M_{\alpha, N}^{(m)}$, for $\alpha = j / 1000$, $j = 0, \ldots, 1000$. The columns show, respectively, $\max_{\alpha\in[0,1]}E_1^{(m)}(\alpha)$, $\max_{\alpha\in[0,1]}E_2^{(m)}(\alpha)$ and $\max_{\alpha\in[0,1]}E_3^{(m)}(\alpha)$, as defined in \eqref{e:errorsEm}, for $m = 1, \ldots, 6$. The experiments corresponding to $v_1$ have been done with $L = 1.6$; those corresponding to $v_2$, with $L = 1.1$, and those corresponding to $v_3$, with $L = 4$. In all cases, $N = 64$.} \label{t:testDam}
\end{table}

In order to improve the results, we have estimated the convergence rate $\mathcal O(m^{-\mu(\alpha)} )$, where, in our case, $\mu(\alpha)$ is given by
\begin{equation}
\label{e:mualpha}
\mu(\alpha) = \lim_{m\to\infty}\log_2\left(\frac{E^{(m)}(\alpha)}{E^{(m+1)}(\alpha)}\right),
\end{equation}

\noindent because, from \eqref{e:int0pim}, $E^{(m+1)}$ requires twice as many points $s_j^{(m)}$ as $E^{(m+1)}$. Note however that, since the accuracy of a computer is finite, we do not take in this formula $m$ tending to infinity, but just an $m$ large enough.
\begin{figure}[htb!]
\centering
 \includegraphics[width=0.48\textwidth,clip=true]{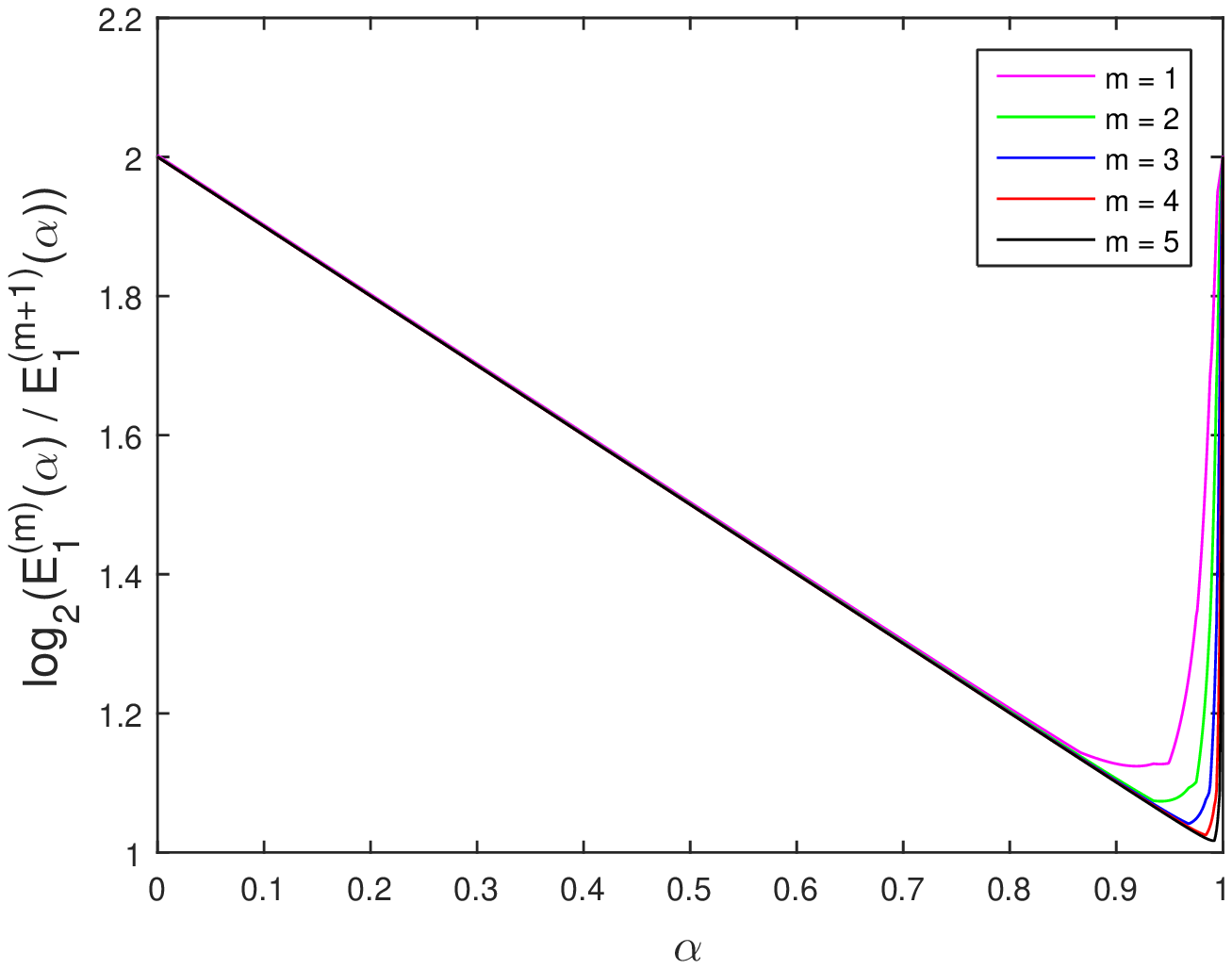}~~~
 \includegraphics[width=0.48\textwidth,clip=true]{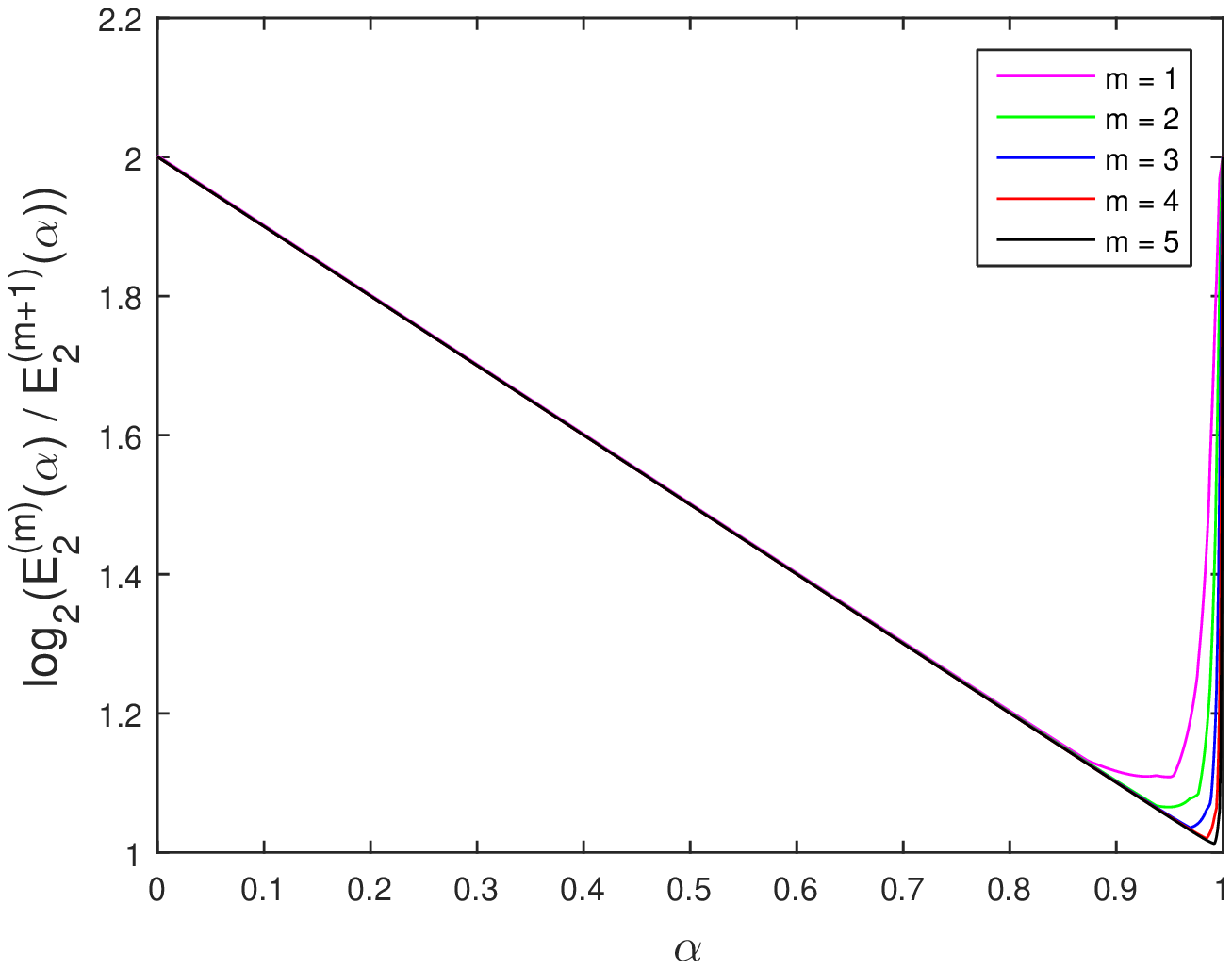}
 \includegraphics[width=0.48\textwidth,clip=true]{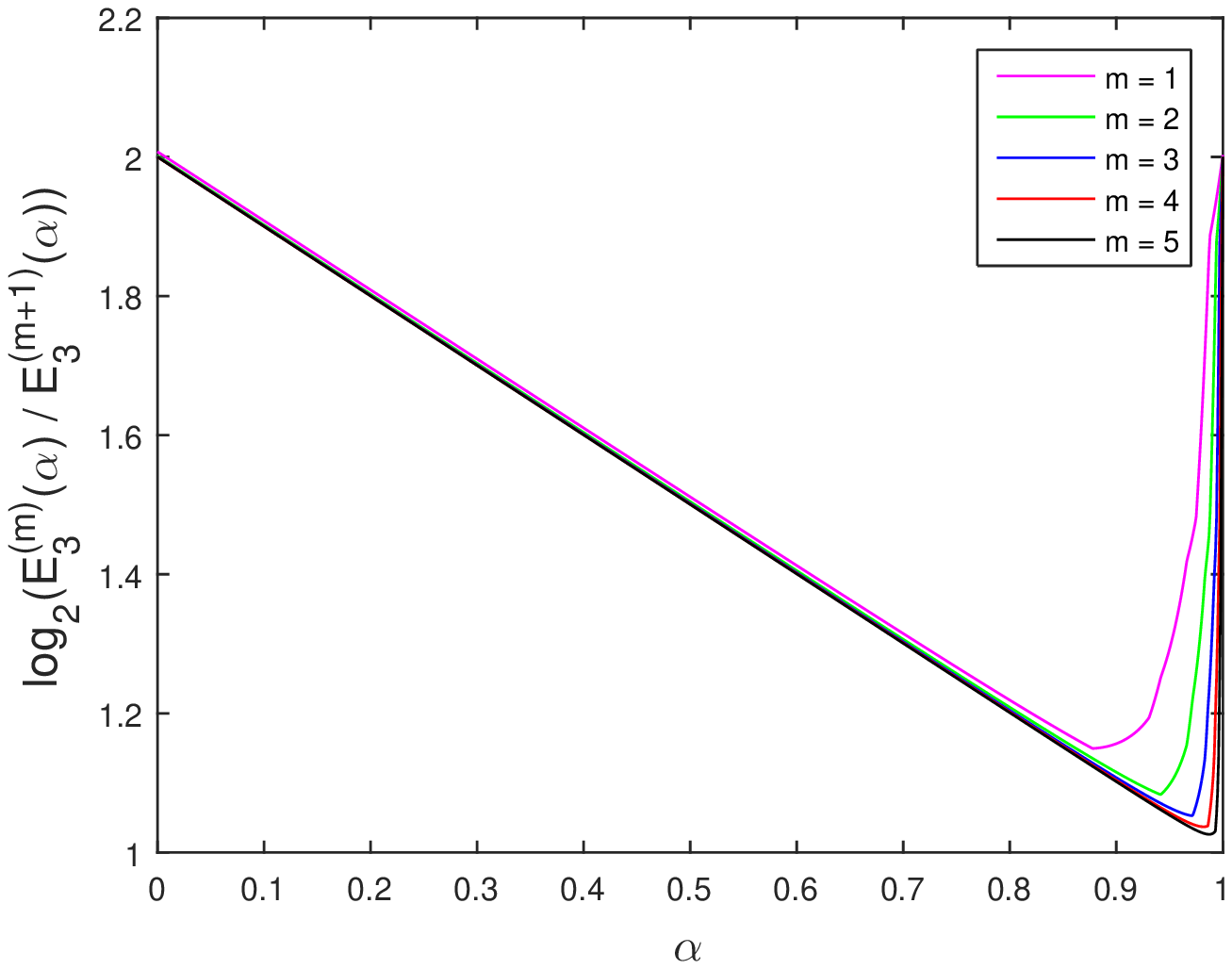}
\caption{Convergence rates of $[\partial_x\mathcal D^\alpha]^{(m)}$ to $\partial_x\mathcal D^\alpha$. Left: $\log_2(E_1^{(m)}(\alpha) / E_1^{(m + 1)}(\alpha))$, corresponding to $v_1(x)$. Right: $\log_2(E_2^{(m)}(\alpha) / E_2^{(m + 1)}(\alpha))$, corresponding to $v_2(x)$. Bottom: $\log_2(E_3^{(m)}(\alpha) / E_3^{(m + 1)}(\alpha))$, corresponding to $v_3(x)$. As $m$ grows, the curves tend to $\mu(\alpha) = 2 - \alpha$, for $\alpha\in[0,1)$.}
\label{f:log2EmEm+1}
\end{figure}

In Figure \ref{t:testDam}, we have plotted the convergence rates corresponding to $v_1(x)$, $v_2(x)$ and $v_3(x)$, which are respectively given by $\log_2(E_1^{(m)}(\alpha) / E_1^{(m + 1)}(\alpha))$ (left), $\log_2(E_2^{(m)}(\alpha) / E_2^{(m + 1)}(\alpha))$ (right), and $\log_2(E_3^{(m)}(\alpha) / E_3^{(m + 1)}(\alpha))$ (bottom). The three graphics are very similar, and give evidence that $\mu(\alpha) = 2-\alpha$. For most $\alpha$, this is clear even for $m = 1$, although, for those $\alpha$ very close to $\alpha = 1$, we can see the curves quickly converge to $\mu(\alpha) = 2 - \alpha$, as $m$ grows. Hence, we claim that
\begin{equation}
\label{e:errorsEmO}
\left\|\left[\partial_x\mathcal D^\alpha\right]^{m} v(x) - \partial_x\mathcal D^\alpha v(x)\right\|_{\infty} =
\mathcal O\left(\frac{1}{m^{2 - \alpha}}\right).
\end{equation}

\noindent Bearing in mind \eqref{e:errorsEmO}, it is immediate to construct more accurate quadrature formulas by applying Richardson extrapolation \cite{richardson} to \eqref{e:int0pim}. More precisely, we define
\begin{equation}
\label{e:extrapolation1}
\left[\partial_x\mathcal D^\alpha\right]^{(m,m+1)} u(s_j) \equiv \frac{2^{2-\alpha}\left[\partial_x\mathcal D^\alpha\right]^{(m+1)}u(s_j) - \left[\partial_x\mathcal D^\alpha\right]^{(m)}u(s_j)}{2^{2-\alpha}-1},
\end{equation}

\noindent together with its associated matrix representation $\mathbf M_{\alpha, N}^{(m,m+1)}$. Note that Figure \ref{t:testDam} suggests that \eqref{e:errorsEmO} is valid for $\alpha\in[0,1)$. On the other hand, when $\alpha = 1$, the experiments show clearly $\mu(1) = 2$, i.e., a second-order convergence. In general, the convergence for the limiting cases $\alpha = 0$ and $\alpha = 1$ is often better. However, since we are interested in $\alpha\in(0,1)$, these non-fractionary cases are not so relevant, so we have applied \eqref{e:extrapolation1} and the extrapolation formulas that will appear in the next pages for all $\alpha\in[0,1]$, obtaining good results even for $\alpha = 0$ and $\alpha = 1$.

As in \eqref{e:errorsEm}, we define $E^{(m,m+1)}(\alpha) \equiv\max_j|[\mathbf M_{\alpha, N}^{(m,m+1)}\cdot\widehat{\mathbf U}]_j - \partial_x\mathcal D^\alpha v(x)|$, and as before, we let $E_i^{(m,m+1)}(\alpha)$ correspond to the error of $v_i(x)$ for each $i=1,2,3$. Table \ref{t:testDamm+1} shows $\max_{\alpha\in[0,1]}E_1^{(m,m+1)}(\alpha)$, $\max_{\alpha\in[0,1]}E_2^{(m,m+1)}(\alpha)$ and $\max_{\alpha\in[0,1]}E_3^{(m,m+1)}(\alpha)$, for $m = 1, \ldots, 5$. Observe that the errors decay faster than in Table \ref{t:testDamm+1}.
\begin{table}[htb!]
\centering
\begin{tabular}{|c|c|c|c|}
\cline{2-4} \multicolumn{1}{c|}{} & $v_1(x)$ & $v_2(x)$ & $v_3(x)$
    \cr
\hline $\mathbf M^{(1,2)}$ &$7.5676\cdot10^{-4}$ &$1.0814\cdot10^{-3}$ &$2.3110\cdot10^{-3}$
    \cr
\hline $\mathbf M^{(2,3)}$ &$1.8889\cdot10^{-4}$ &$2.7010\cdot10^{-4}$ &$5.7575\cdot10^{-4}$
    \cr
\hline $\mathbf M^{(3,4)}$ &$4.7205\cdot10^{-5}$ &$6.7510\cdot10^{-5}$ &$1.4381\cdot10^{-4}$
    \cr
\hline $\mathbf M^{(4,5)}$ &$1.1800\cdot10^{-5}$ &$1.6877\cdot10^{-5}$ &$3.5945\cdot10^{-5}$
    \cr
\hline $\mathbf M^{(5,6)}$ &$2.9499\cdot10^{-6}$ &$4.2191\cdot10^{-6}$ &$8.9858\cdot10^{-6}$
    \cr
\hline
\end{tabular}
\caption{Continuation of Table \ref{t:testDam}. Errors obtained after computing the fractional derivatives via $\mathbf M_{\alpha, N}^{(m,m+1)}$, for $m = 1, \ldots, 5$. All the other details are identical.} \label{t:testDamm+1}
\end{table}

\begin{figure}[htb!]
\centering
 \includegraphics[width=0.48\textwidth,clip=true]{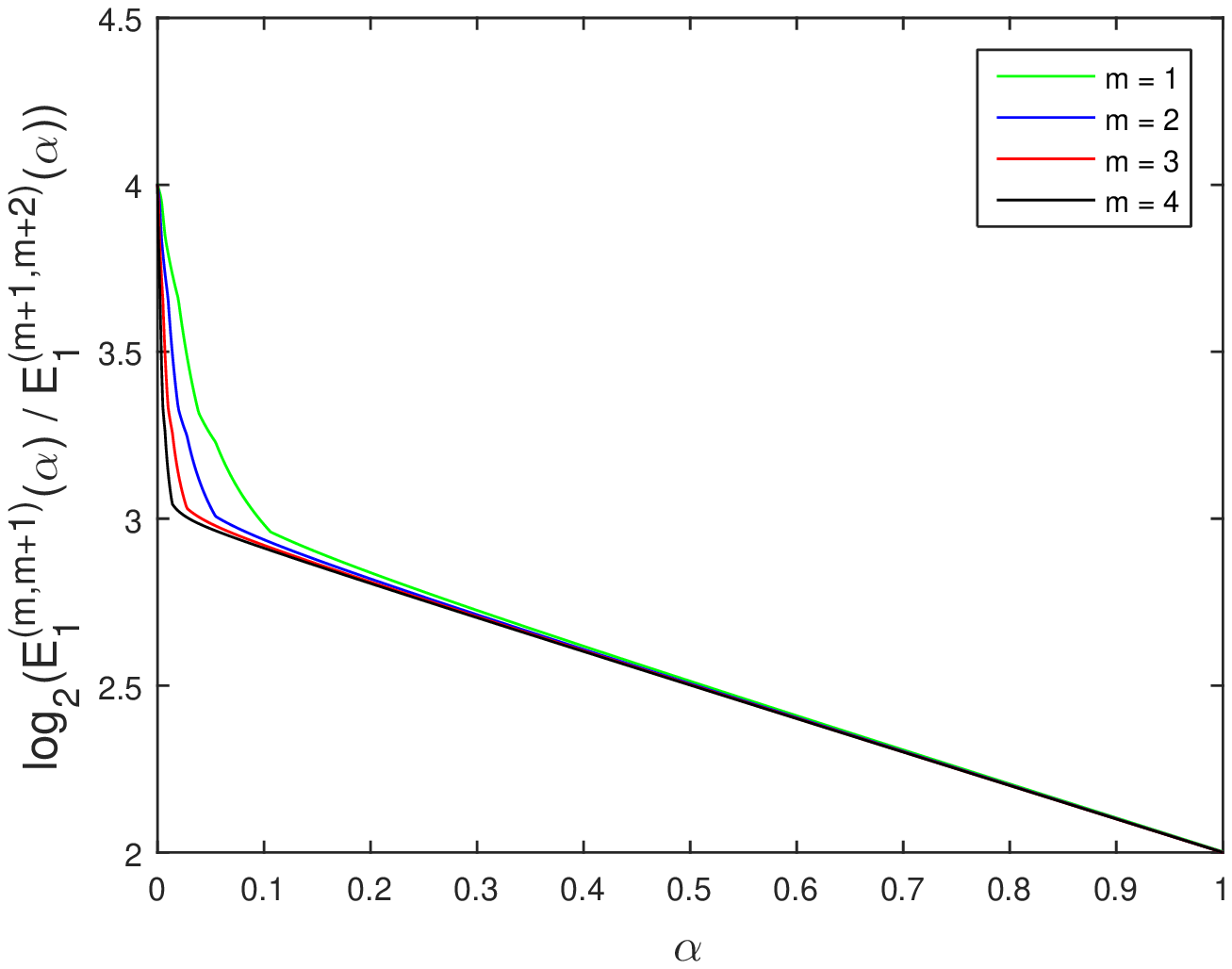}~~~
 \includegraphics[width=0.48\textwidth,clip=true]{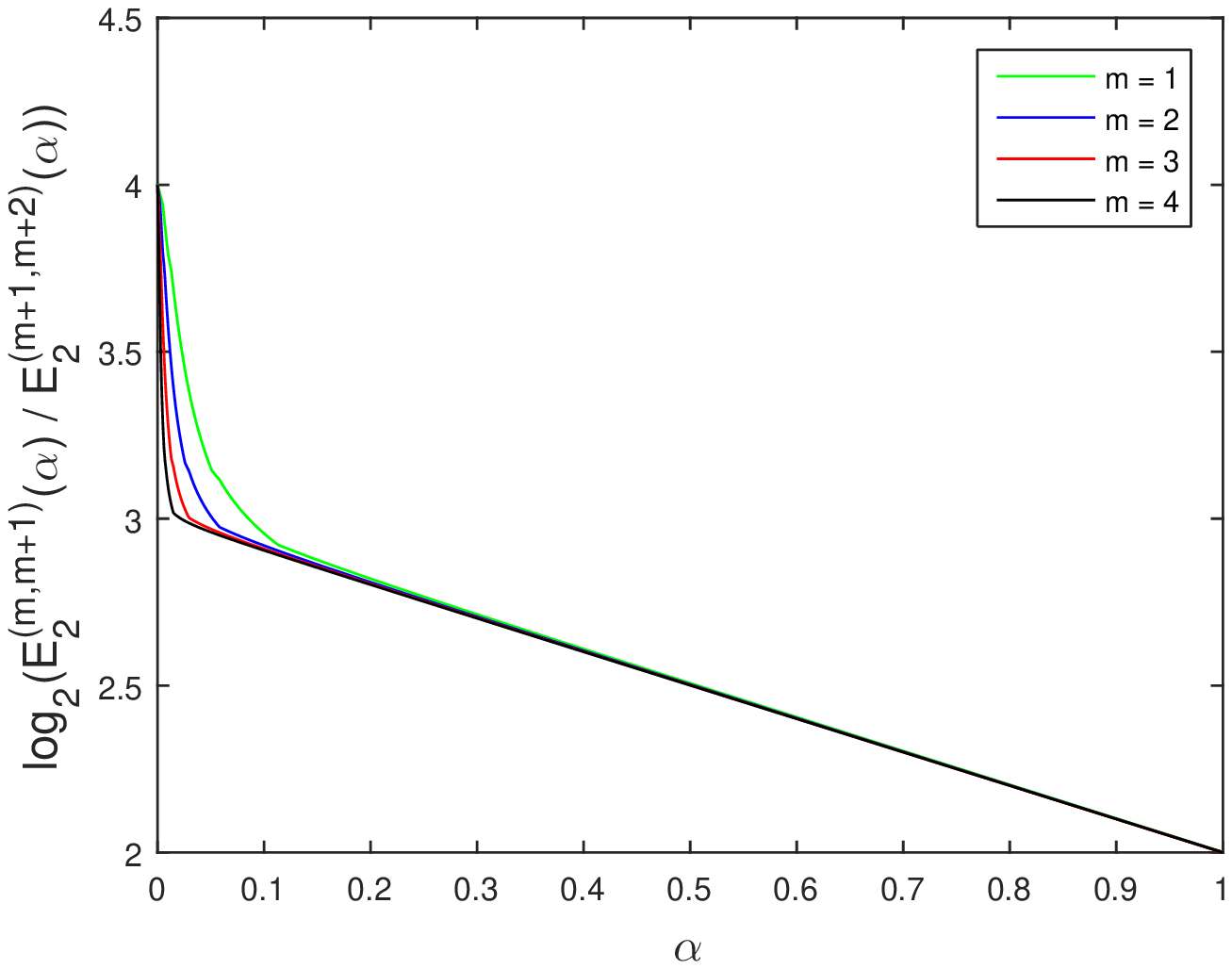}
 \includegraphics[width=0.48\textwidth,clip=true]{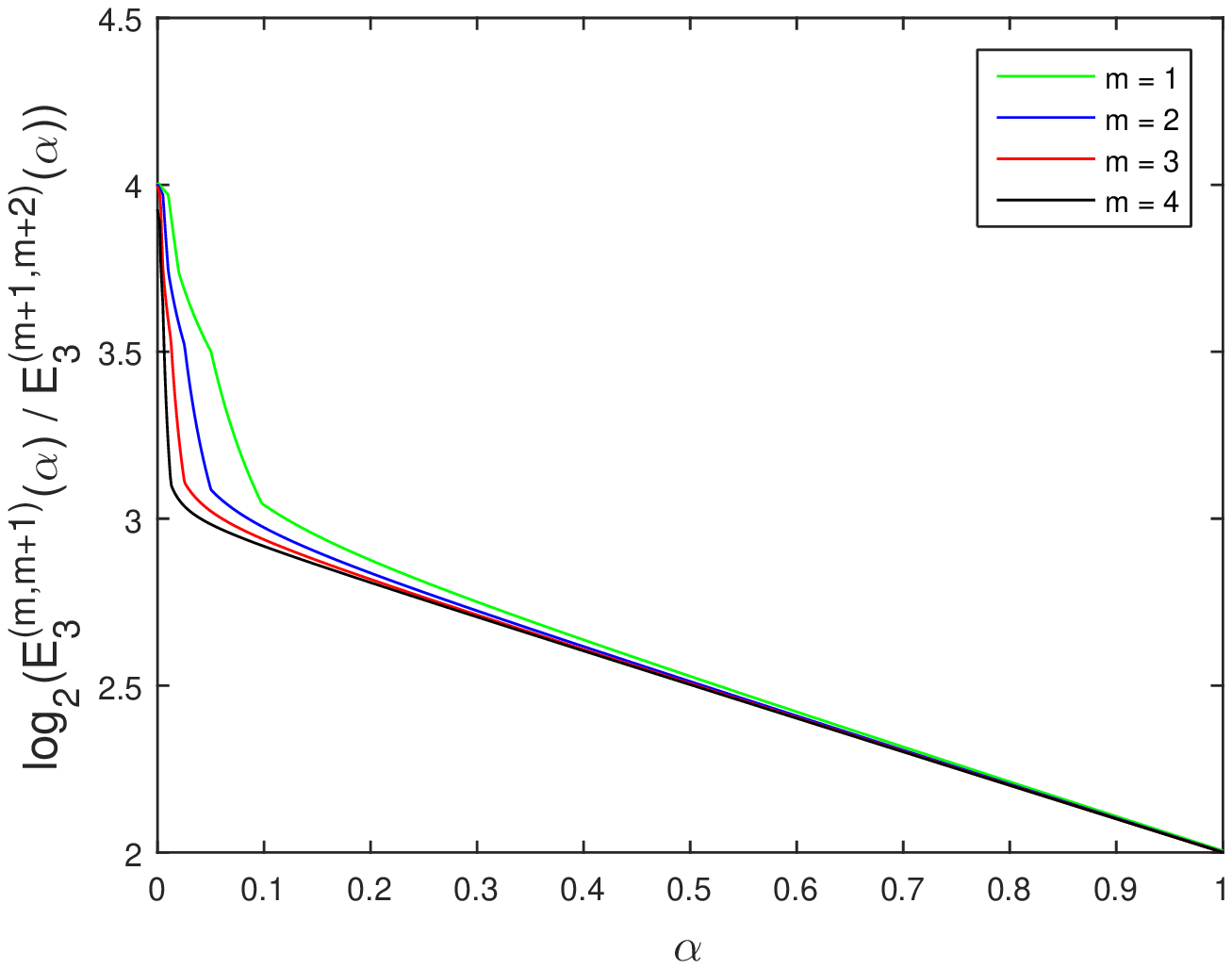}
\caption{Convergence rates for \eqref{e:extrapolation1}. Left: $\log_2(E_1^{(m,m+1)}(\alpha) / E_1^{(m+1,m+2)}(\alpha))$, corresponding to $v_1(x)$. Right: $\log_2(E_2^{(m,m+1)}(\alpha) / E_2^{(m+1,m+2)}(\alpha))$, corresponding to $v_2(x)$. Bottom: $\log_2(E_3^{(m,m+1)}(\alpha) / E_3^{(m+1,m+2)}(\alpha))$, corresponding to $v_3(x)$. As $m$ grows, the curves tend to $\mu(\alpha) = 3 - \alpha$, for $\alpha\in(0,1]$.}
\label{f:log2EmEm+1m+2}
\end{figure}

Figure \ref{f:log2EmEm+1m+2} depicts the convergence rates corresponding to Table \ref{t:testDamm+1}, given by $\log_2(E_1^{(m,m+1)}(\alpha) / E_1^{(m+1,m+2)}(\alpha))$ (left), $\log_2(E_2^{(m,m+1)}(\alpha) / E_2^{(m+1,m+2)}(\alpha))$ (right) and $\log_2(E_3^{(m,m+1)}(\alpha) / E_3^{(m+1,m+2)}(\alpha))$ (bottom), which correspond respectively to $v_1(x)$, $v_2(x)$ and $v_3(x)$. The three graphics are very similar, and give evidence that $\mu(\alpha) = 3-\alpha$. For most $\alpha$, this is clear even for $m = 1$, although, for those $\alpha$ very close to $\alpha = 0$, we can see the curves quickly converge to $\mu(\alpha) = 3 - \alpha$ as $m$ grows. Thus, there is a complete parallelism with Figure \ref{f:log2EmEm+1}.

In view of the previous results, we claim that
\begin{equation}
\label{e:errorsEmm+1O}
\left\|\left[\partial_x\mathcal D^\alpha\right]^{m,m+1} v(x) - \partial_x\mathcal D^\alpha v(x)\right\|_{\infty} =
\mathcal O\left(\frac{1}{m^{3 - \alpha}}\right);
\end{equation}

\noindent hence,
\begin{equation*}
\left\|\left[\partial_x\mathcal D^\alpha\right]^{m} v(x) - \partial_x\mathcal D^\alpha v(x)\right\|_{\infty} =
\frac{c_1(\alpha)}{m^{2 - \alpha}} + \frac{c_2(\alpha)}{m^{3 - \alpha}} + \ldots,
\end{equation*}

\noindent for some $c_1(\alpha)$, $c_2(\alpha)$ bounded on $\alpha\in(0,1)$. Bearing in mind \eqref{e:errorsEmm+1O}, we have applied Richardson extrapolation to \eqref{e:extrapolation1}, to define
\begin{equation}
\label{e:extrapolation2}
\left[\partial_x\mathcal D^\alpha\right]^{(m,m+1,m+2)} u(s_j) \equiv \frac{2^{3-\alpha}\left[\partial_x\mathcal D^\alpha\right]^{(m+1,m+2)}u(s_j) - \left[\partial_x\mathcal D^\alpha\right]^{(m,m+1)}u(s_j)}{2^{3-\alpha}-1},
\end{equation}

\noindent together with its associated matrix representation $\mathbf M_{\alpha, N}^{(m,m+1,m+2)}$.
\begin{table}[htb!]
\centering
\begin{tabular}{|c|c|c|c|}
\cline{2-4} \multicolumn{1}{c|}{} & $v_1(x)$ & $v_2(x)$ & $v_3(x)$
    \cr
\hline $\mathbf M^{(1,2,3)}$ &$5.5822\cdot10^{-7}$ &$6.0120\cdot10^{-7}$ &$3.4331\cdot10^{-6}$
    \cr
\hline $\mathbf M^{(2,3,4)}$ &$6.0166\cdot10^{-8}$ &$6.6696\cdot10^{-8}$ &$3.6703\cdot10^{-7}$
    \cr
\hline $\mathbf M^{(3,4,5)}$ &$6.7540\cdot10^{-9}$ &$7.5771\cdot10^{-9}$ &$4.0560\cdot10^{-8}$
    \cr
\hline $\mathbf M^{(4,5,6)}$ &$7.6782\cdot10^{-10}$ &$8.6835\cdot10^{-10}$ &$4.5490\cdot10^{-9}$
    \cr
\hline
\end{tabular}
\caption{Continuation of Tables \ref{t:testDam} and \ref{t:testDamm+1}. Errors obtained after computing the fractional derivatives via $\mathbf M_{\alpha, N}^{(m,m+1,m+2)}$, for $m = 1, \ldots, 4$. All the other details are identical.}\label{t:testDamm+1m+2}
\end{table}

Table \ref{t:testDamm+1m+2} shows $\max_{\alpha\in[0,1]}E_1^{(m,m+1,m+2)}(\alpha)$, $\max_{\alpha\in[0,1]}E_2^{(m,m+1,m+2)}(\alpha)$ and $\max_{\alpha\in[0,1]}E_3^{(m,m+1,m+2)}(\alpha)$, for $m = 1, \ldots, 4$; where $E^{(m,m+1,m+2)}(\alpha)$ is the generalization of \eqref{e:errorsEm} corresponding to $\mathbf M_{\alpha, N}^{(m,m+1,m+2)}$, etc. The errors decay now much faster than in Tables \ref{t:testDam} and \ref{t:testDamm+1}. Moreover, in view of \eqref{e:errorsEmO} and \eqref{e:errorsEmm+1O}, we may expect that
\begin{equation}
\label{e:errorsEmm+1m+2O}
\left\|\left[\partial_x\mathcal D^\alpha\right]^{m,m+1,m+2} v(x) - \partial_x\mathcal D^\alpha v(x)\right\|_{\infty} = \mathcal O\left(\frac{1}{m^{4 - \alpha}}\right).
\end{equation}

\begin{figure}[htb!]
\centering
 \includegraphics[width=0.48\textwidth,clip=true]{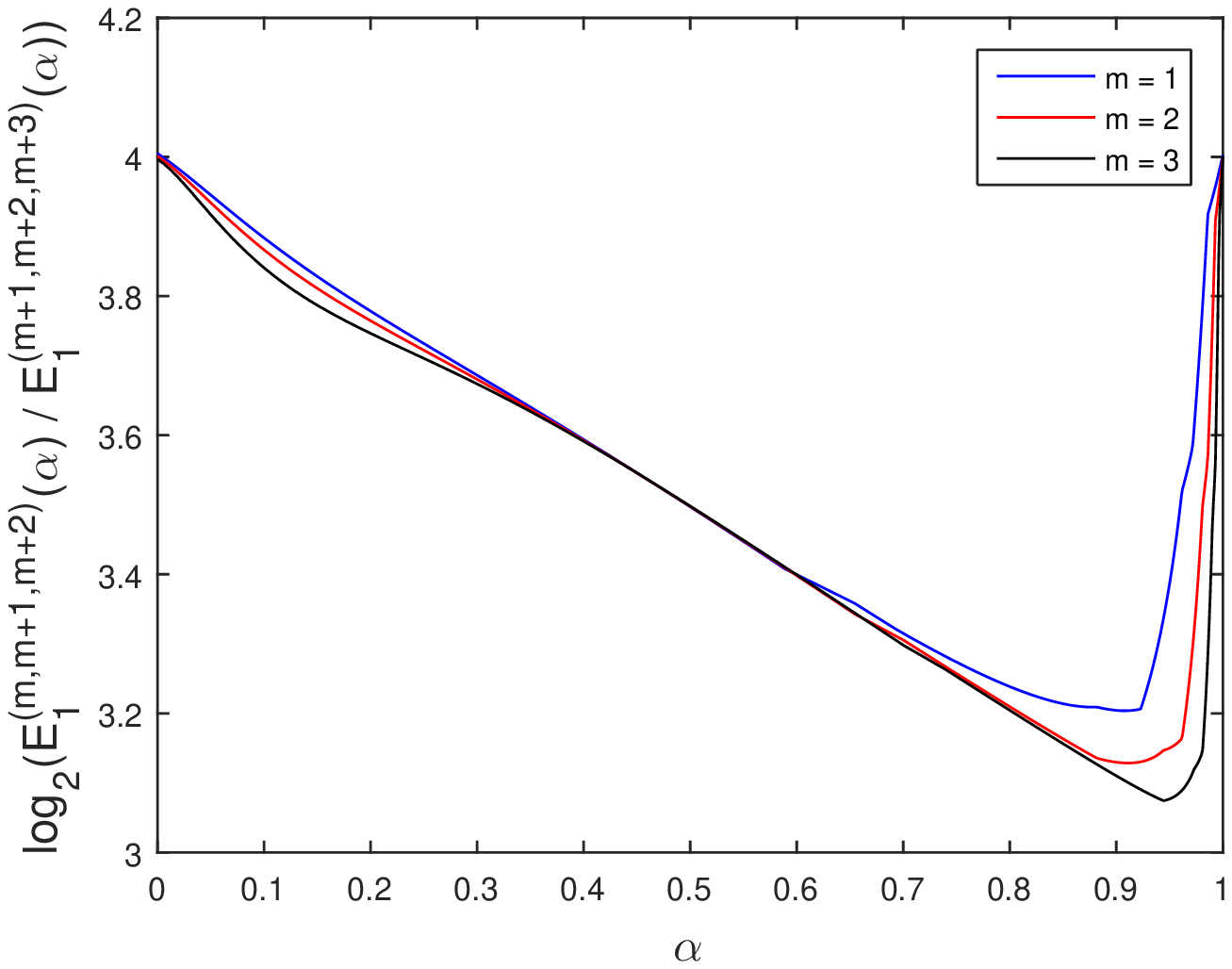}~~~
 \includegraphics[width=0.48\textwidth,clip=true]{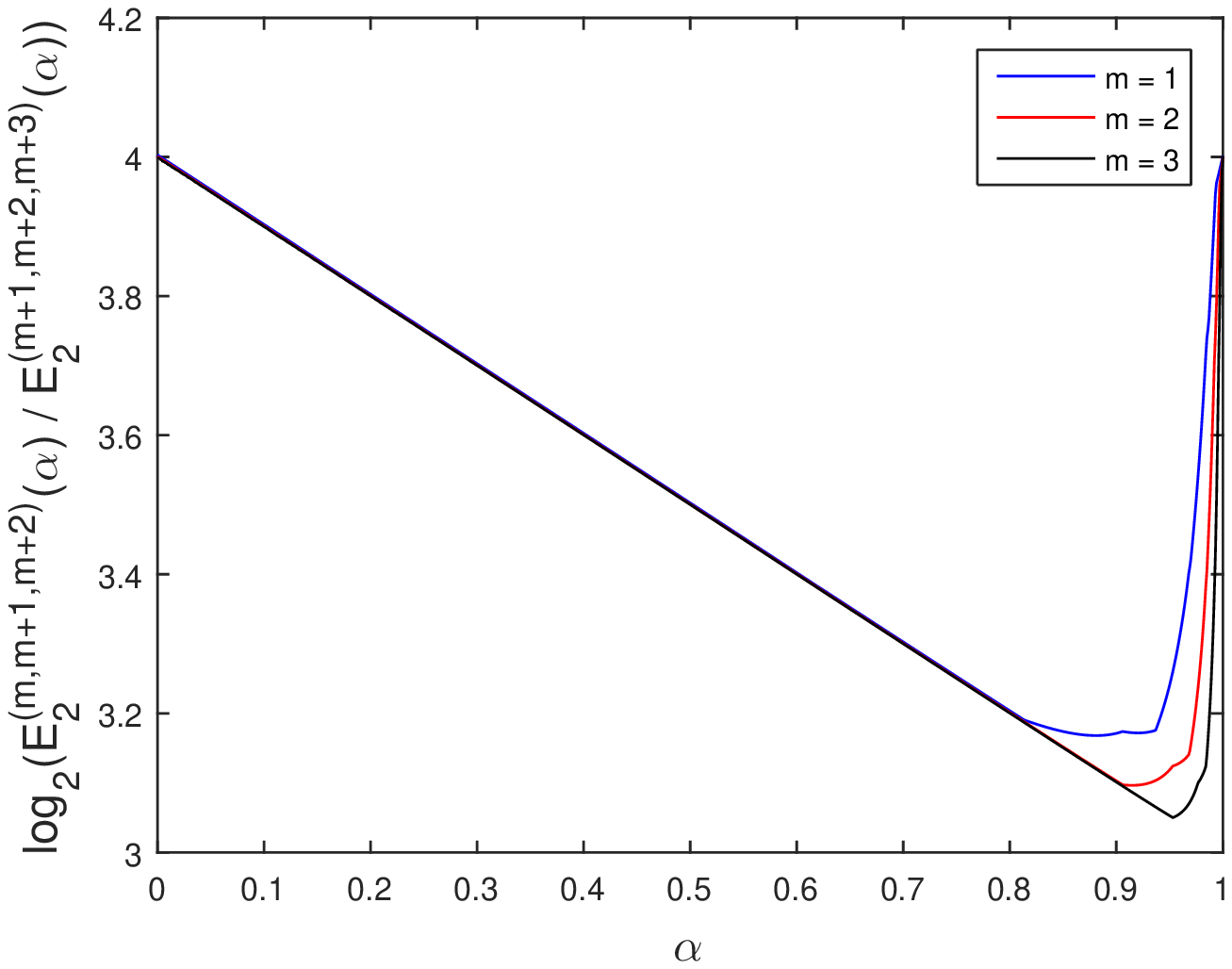}
 \includegraphics[width=0.48\textwidth,clip=true]{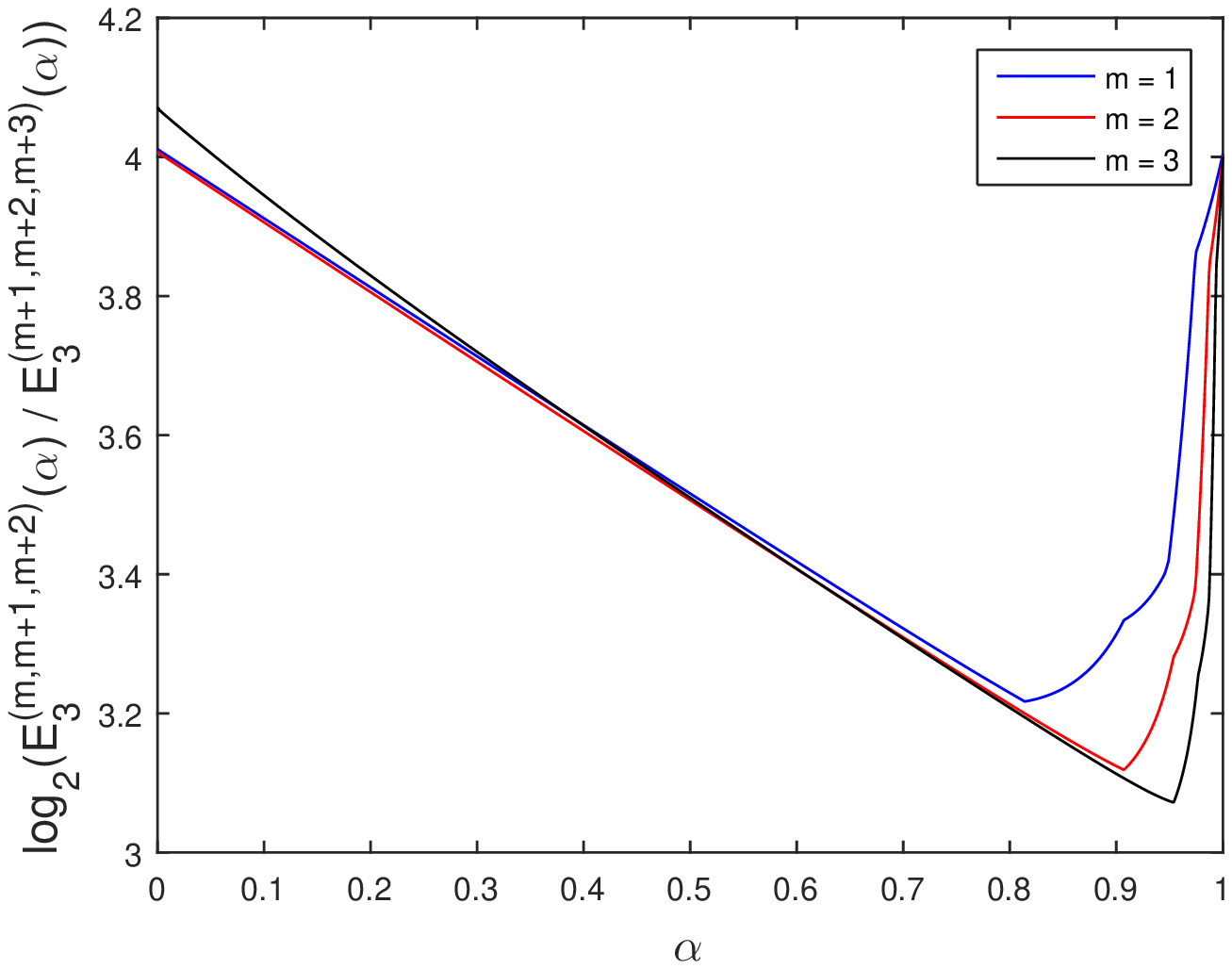}
\caption{Convergence rates for \eqref{e:extrapolation2}. Left: $\log_2(E_1^{(m,m+1,m+2)}(\alpha) / E_1^{(m+1,m+2,m+3)}(\alpha))$, corresponding to $v_1(x)$. Right: $\log_2(E_2^{(m,m+1,m+2)}(\alpha) / E_2^{(m+1,m+2,m+3)}(\alpha))$, corresponding to $v_2(x)$. Bottom: $\log_2(E_3^{(m,m+1,m+2)}(\alpha) / E_3^{(m+1,m+2,m+3)}(\alpha))$, corresponding to $v_3(x)$. The curves can be roughly approximated by $\mu(\alpha) = 4 - \alpha$, for $\alpha\in[0, 1)$.}
\label{f:log2EmEm+1m+2m+3}
\end{figure}

\noindent This is confirmed by Figure \ref{f:log2EmEm+1m+2m+3}, which depicts the convergence rates corresponding to Table \ref{t:testDamm+1m+2}, given by $\log_2(E_1^{(m,m+1,m+2)}(\alpha) / E_1^{(m+1,m+2,m+3)}(\alpha))$ for $v_1(x)$ (left); by $\log_2(E_2^{(m,m+1,m+2)}(\alpha) / E_2^{(m+1,m+2,m+3)}(\alpha))$ for $v_2(x)$ (right); and by $\log_2(E_3^{(m,m+1,m+2)}(\alpha) / E_3^{(m+1,m+2,m+3)}(\alpha))$ for $v_3(x)$ (bottom). Since the errors are now much smaller, the graphics are not so sharp as those in Figures \ref{f:log2EmEm+1} and \ref{f:log2EmEm+1m+2}. Nevertheless, they still give acceptable evidence that $\mu(\alpha) = 4 - \alpha$.

\noindent We can go further in the extrapolation process. Indeed, all the previous arguments strongly suggest that
\begin{equation}
\label{e:assympt}
\begin{split}
\|[\partial_x\mathcal D^\alpha]^{m} v(x) & - \partial_x\mathcal D^\alpha v(x)\|_{\infty}
    \cr
& =
\frac{c_1(\alpha)}{m^{2 - \alpha}} + \frac{c_2(\alpha)}{m^{3 - \alpha}} + \frac{c_3(\alpha)}{m^{4 - \alpha}} + \frac{c_4(\alpha)}{m^{5 - \alpha}} + \frac{c_5(\alpha)}{m^{6 - \alpha}} + \ldots,
\end{split}
\end{equation}

\noindent for some $c_1(\alpha)$, $c_2(\alpha)$, $c_3(\alpha)\ldots$ bounded on $\alpha\in(0,1)$. Therefore, we define, for a given integer $n\ge 3$,
\begin{equation}
\label{e:extrapolationn}
\begin{split}
[\partial_x\mathcal D^\alpha] & ^{(m,m+1,\ldots,m + n - 1, m + n)} u(s_j)
    \cr
& \equiv \frac{2^{n + 1 -\alpha}\left[\partial_x\mathcal D^\alpha\right]^{(m+1,\ldots,m+n)}u(s_j) - \left[\partial_x\mathcal D^\alpha\right]^{(m,\ldots,m + n - 1)}u(s_j)}{2^{n + 1 -\alpha}-1},
\end{split}
\end{equation}

\noindent together with its associated matrix representation $\mathbf M_{\alpha, N}^{(m,m+1,\ldots,m + n - 1, m + n)}$. In Table \ref{t:testDamhigher}, we have used higher-order approximations given by \eqref{e:extrapolationn} to approximate $\partial_x\mathcal D^\alpha$. The very high accuracy achieved with, for instance, $\mathbf M^{(1,2,3,4,5)}$ and $\mathbf M^{(1,2,3,4,5,6)}$, confirms the adequacy of \eqref{e:assympt}.

\begin{table}[htb!]
\centering
\begin{tabular}{|c|c|c|c|}
\cline{2-4} \multicolumn{1}{c|}{} & $v_1(x)$ & $v_2(x)$ & $v_3(x)$
    \cr
\hline $\mathbf M^{(1,2,3,4)}$ &$2.8383\cdot10^{-8}$ &$2.3888\cdot10^{-8}$ &$1.9085\cdot10^{-7}$
    \cr
\hline $\mathbf M^{(2,3,4,5)}$ &$1.7682\cdot10^{-9}$ &$1.4907\cdot10^{-9}$ &$1.1862\cdot10^{-8}$
    \cr
\hline $\mathbf M^{(3,4,5,6)}$ &$1.1053\cdot10^{-10}$ &$9.3191\cdot10^{-11}$ &$7.4030\cdot10^{-10}$
    \cr
\hline $\mathbf M^{(1,2,3,4,5)}$ &$7.7385\cdot10^{-12}$ &$3.5603\cdot10^{-12}$ &$7.6762\cdot10^{-11}$
    \cr
\hline $\mathbf M^{(2,3,4,5,6)}$ &$6.6334\cdot10^{-13}$ &$3.8475\cdot10^{-13}$ &$8.7386\cdot10^{-12}$
    \cr
\hline $\mathbf M^{(1,2,3,4,5,6)}$ &$5.3305\cdot10^{-13}$ &$3.8948\cdot10^{-13}$ &$7.7834\cdot10^{-12}$
    \cr
\hline
\end{tabular}
\caption{Continuation of Tables \ref{t:testDam}, \ref{t:testDamm+1} and \ref{t:testDamm+1m+2}. Errors obtained after computing the fractional derivatives via higher-order approximations. All the other details are identical.}\label{t:testDamhigher}
\end{table}

\begin{figure}[htb!]
\centering
\includegraphics[width=0.48\textwidth,clip=true]{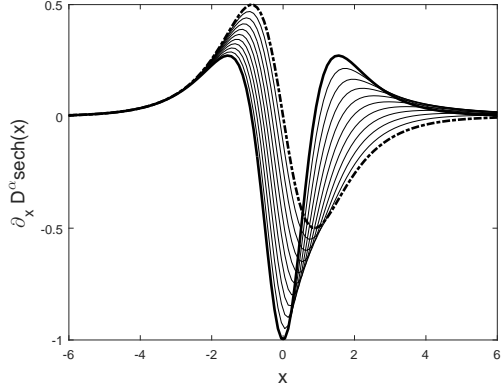}
\caption{
Numerical computation of $\partial_x\mathcal D^\alpha\sech(x)$, for $\alpha = 0, 0.1, \ldots, 1$. The thick dashed-dotted line represents the computation for $\alpha = 0$ and the thick solid line represents the computation for $\alpha = 1$; being the thin solid lines the transitions between these two.
}
\label{f:frac01}
\end{figure}

It is interesting to show graphically how fractional derivatives look like for different values of $\alpha$. In Figure \ref{f:frac01}, we have plotted $\left[\partial_x\mathcal D^{\alpha}\right]^{(1,2,3,4,5,6)}\sech(x)$, for $\alpha = 0, 0.1, \ldots, 1$, taking $N= 128$, $L = 3.9$. Note that, unlike in the previous examples, we do not know the explicit expression of $\partial_x\mathcal D^\alpha \sech(x)$, except in the limiting cases $\alpha = 0$ and $\alpha = 1$, which correspond to $v'(x) = -\tanh(x)/\cosh(x)$ and $v''(x) = (2\tanh^2(x) - 1) / \cosh(x)$, respectively. However, since $\|\left[\partial_x\mathcal D^{\alpha = 0}\right]^{(1,2,3,4,5,6)}\sech(x) - \sech'(x)\|_\infty = 2.3873\cdot10^{-13}$ and $\|\left[\partial_x\mathcal D^{\alpha = 1}\right]^{(1,2,3,4,5,6)}\sech(x) - \sech''(x)\|_\infty = 7.2609 \cdot10^{-14}$, we can expect a similar accuracy also for $\alpha\in(0,1)$.

Observe the transition between the thick dashed-dotted line ($\alpha = 0$, i.e., $\sech'(x)$) and the thick solid line ($\alpha = 1$, i.e., $\sech''(x)$). Let us recall here the comments on the related Figure 5 from \cite{hale2012}: it is rare to see figures like Figure \ref{f:frac01}, in which numerically evaluated fractional derivatives or integrals are plotted. Indeed, the authors in \cite{hale2012} claimed that they did not known of such figures in the literature. Therefore, Figure \ref{f:frac01} seems to be the first one where the numerically evaluated fractional derivative of a function defined over the whole $\mathbb R$ has been ever plotted.

\subsection{Remarks on the minimum regularity required}

\label{s:regularity}

It is well-known that pseudo-spectral methods are most convenient when applied to regular functions. In fact, along this paper, we are assuming that all the functions are $\mathcal C^3$, and satisfying \eqref{e:decay1} and \eqref{e:decay2}. However, it is straightforward to check that our method can still be useful for functions with much less regularity. In order to illustrate this, we have considered four functions with similar structure, and decaying quadratically at infinity, but with growing regularity:
\begin{equation}
\label{e:v4567}
\begin{split}
v_4(x) & = \frac{x|x|}{1 + x^4}, \qquad v_5(x) = \frac{|x|^3}{1 + |x|^5},
    \cr
v_6(x) & = \frac{x^3|x|}{1 + x^6}, \qquad v_7(x) = \frac{|x|^5}{1 + |x|^7}.
\end{split}
\end{equation}

\noindent More precisely, $v_4(x)$ is $\mathcal C^1$, but not $\mathcal C^2$; $v_5(x)$ is $\mathcal C^2$, but not $\mathcal C^3$; $v_6(x)$ is $\mathcal C^3$, but not $\mathcal C^4$; and $v_7(x)$ is $\mathcal C^4$, but not $\mathcal C^5$.

Unlike in the previous examples with $v_1(x)$, $v_2(x)$ and $v_3(x)$, we do not know the explicit expression of $\partial_x\mathcal D^\alpha$ applied to \eqref{e:v4567}. Therefore, as in the example with $\sech(x)$ in Figure \ref{f:frac01}, we have to content ourselves with testing our method for the limiting cases $\alpha = 0$ and $\alpha = 1$, which correspond respectively to the first and second derivatives:
\begin{equation*}
\begin{split}
v_4'(x) & = -\frac{2|x|^5 - 2|x|}{(1 + x^4)^2}, \qquad v_4''(x) = \frac{6x^8 - 24x^4 + 2}{(1 + x^4)^3}\sgn(x),
    \cr
v_5'(x) & = -\frac{2x^7 - 3x|x|}{(1 + |x|^5)^2}, \qquad v_5''(x) = \frac{6|x|^{11} - 38x^6 + 6|x|}{(1 + |x|^5)^3},
    \cr
v_6'(x) & = -\frac{2|x|^9 - 4|x|^3}{(1 + x^6)^2}, \qquad v_6''(x) = \frac{6x^{13} - 54x^7 + 12x}{(1 + x^6)^3}|x|,
    \cr
v_7'(x) & = -\frac{2x^{11} - 5x^3|x|}{(1 + |x|^7)^2}, \qquad v_7''(x) = \frac{6|x|^{17} - 72x^{10} + 20|x|^3}{(1 + |x|^7)^3}.
\end{split}
\end{equation*}

\noindent Table \ref{t:testC4567} shows the errors in discrete $L^\infty$-norm between the exact and the approximated values of $\partial_x\mathcal D^\alpha v_4(x)$, $\partial_x\mathcal D^\alpha v_5(x)$, $\partial_x\mathcal D^\alpha v_6(x)$, and $\partial_x\mathcal D^\alpha v_7(x)$, for $\alpha = 0$ and $\alpha = 1$, and for six different-order approximations of $\partial_x\mathcal D^\alpha$. The experiments have been done respectively with $L = 0.1$, $L = 0.14$, $L = 0.2$ and $L = 0.28$. In all cases, $N = 256$.
\begin{table}[htb!]
    \centering
\begin{tabular}{|c|c|c||c|c|c|}
    \cline{2-5}
    \multicolumn{1}{c|}{}& \multicolumn{2}{c||}{$v_4(x)$} & \multicolumn{2}{c|}{$v_5(x)$}
\cr
    \cline{2-5} \multicolumn{1}{c|}{} & $\alpha=0$ & $\alpha = 1$ & $\alpha = 0$ & $\alpha = 1$
\cr
    \hline $\mathbf M^{(1)}$ &$1.7691\cdot10^{-3}$ &$1.0451\cdot10^{-1}$ &$1.6052\cdot10^{-3}$ &$8.1360\cdot10^{-3}$
\cr
    \hline $\mathbf M^{(1,2)}$ &$1.8867\cdot10^{-4}$ &$8.6964\cdot10^{-2}$ &$4.1015\cdot10^{-6}$ &$4.0818\cdot10^{-3}$
\cr
    \hline $\mathbf M^{(1,2,3)}$ &$1.8823\cdot10^{-4}$ &$8.4886\cdot10^{-2}$ &$4.4278\cdot10^{-7}$ &$5.7049\cdot10^{-4}$
\cr
    \hline $\mathbf M^{(1,2,3,4)}$ &$1.8827\cdot10^{-4}$ &$8.4396\cdot10^{-2}$ &$4.4260\cdot10^{-7}$ &$5.6994\cdot10^{-4}$
\cr
    \hline $\mathbf M^{(1,2,3,4,5)}$ &$1.8827\cdot10^{-4}$ &$8.4430\cdot10^{-2}$ &$4.4261\cdot10^{-7}$ &$5.6998\cdot10^{-4}$
\cr
    \hline $\mathbf M^{(1,2,3,4,5,6)}$ &$1.8827\cdot10^{-4}$ &$8.4430\cdot10^{-2}$ &$4.4261\cdot10^{-7}$ &$5.6998\cdot10^{-4}$
\cr
    \hline
\multicolumn{5}{c}{}
\cr
\cline{2-5}\multicolumn{1}{c|}{}& \multicolumn{2}{c||}{$v_6(x)$} & \multicolumn{2}{c|}{$v_7(x)$}
    \cr
\cline{2-5} \multicolumn{1}{c|}{} & $\alpha=0$ & $\alpha = 1$ & $\alpha = 0$ & $\alpha = 1$
    \cr
\hline $\mathbf M^{(1)}$ &$1.2645\cdot10^{-3}$ &$8.3822\cdot10^{-3}$ &$9.9656\cdot10^{-4}$ &$8.8407\cdot10^{-3}$
    \cr
\hline $\mathbf M^{(1,2)}$ &$2.5273\cdot10^{-6}$ &$4.2022\cdot10^{-3}$ &$1.3500\cdot10^{-6}$ &$4.4287\cdot10^{-3}$
    \cr
\hline $\mathbf M^{(1,2,3)}$ &$1.8058\cdot10^{-7}$ &$6.4987\cdot10^{-6}$ &$9.6944\cdot10^{-8}$ &$4.9026\cdot10^{-6}$
    \cr
\hline $\mathbf M^{(1,2,3,4)}$ &$1.5933\cdot10^{-9}$ &$1.0984\cdot10^{-6}$ &$5.4074\cdot10^{-11}$ &$3.5168\cdot10^{-7}$
    \cr
\hline $\mathbf M^{(1,2,3,4,5)}$ &$1.5932\cdot10^{-9}$ &$1.0954\cdot10^{-6}$ &$1.1625\cdot10^{-10}$ &$1.0714\cdot10^{-8}$
    \cr
\hline $\mathbf M^{(1,2,3,4,5,6)}$ &$1.5933\cdot10^{-9}$ &$1.0954\cdot10^{-6}$ &$2.2295\cdot10^{-11}$ &$1.0712\cdot10^{-8}$
    \cr
\hline\end{tabular}
\caption{Given $v_4(x) = x|x|(1 + x^4)^{-1}$, $v_5(x) = |x|^3(1 + |x|^5)^{-1}$, $v_6(x) = x^3|x|(1 + x^6)^{-1}$ and $v_7(x) = |x|^5(1 + |x|^7)^{-1}$; errors in discrete $L^\infty$-norm in $s\in[0,\pi]$ between the exact and the approximated values of $\partial_x\mathcal D^\alpha v_4(x)$, $\partial_x\mathcal D^\alpha v_5(x)$, $\partial_x\mathcal D^\alpha v_6(x)$, and $\partial_x\mathcal D^\alpha v_7(x)$, for $\alpha = 0$ and $\alpha = 1$, and for six different-order approximations of $\partial_x\mathcal D^\alpha$. The experiments have been done respectively with $L = 0.1$, $L = 0.14$, $L = 0.2$ and $L = 0.28$. In all cases, $N = 256$.} \label{t:testC4567}
\end{table}

In view of the results, we can draw several important conclusions. On the one hand, as expected, they are not as good as those obtained previously for all $\alpha\in[0,1]$, for the $\mathcal C^\infty$ functions $v_1(x)$, $v_2(x)$ and $v_3(x)$, with just $N = 64$. On the other hand, the accuracy with which one can compute numerically $\partial_x\mathcal D^\alpha v(x)$ for a function $v(x)$ increases dramatically with the regularity of $v(x)$. At this point, we find especially remarkable what happens with $v_4(x)$ when $\alpha = 1$. Indeed, in spite of being approximating a function that is not even continuous,
\begin{equation*}
\lim_{x\to0^-}v_4''(x) = -2, \qquad \lim_{x\to0^+}v_4''(x) = 2,
\end{equation*}

\noindent the results, even if poor, are still coherent.

Let us finish this section by making some comments on the boundedness of the functions. In fact, all the functions we are interested in are bounded, and boundedness is implicitly assumed all the time. However, if a regular function $v(x)$ is unbounded, the change of variable \eqref{e:changeofvariable} will introduce an artificial singularity in $u(s)$ at $s = 0$ and/or $s = \pi$. In those cases, the method is of limited application, but can still be of some utility, at least for functions with \emph{slow} growth. Let us consider for instance
\begin{equation*}
v_8(x) = \log(1 + x^2),
\end{equation*}

\noindent with logarithmic growth, and such that its fractional derivative is explicitly known for all $\alpha$:
\begin{equation}
\label{e:dav8}
\begin{split}
\partial_x \mathcal D^\alpha v_8(x) = & \frac{2\pi\alpha\csc(\alpha\pi)}{\Gamma(1-\alpha)}(1 + x^2)^{-1-\alpha/2}
    \cr
& \cdot\Big[\sin\Big(\frac{\alpha \pi}{2} + \alpha\arctan(x)\Big) + x \cos\Big(\frac{\alpha \pi}{2} + \alpha\arctan(x)\Big)\Big];
\end{split}
\end{equation}

\noindent As with $v_1(x)$, $v_2(x)$ and $v_3(x)$, we have taken $1001$ equally-spaced values of $\alpha\in[0,1]$, i.e., $\alpha = j / 1000$, $j = 0, \ldots, 1000$. When $\alpha = 0$, and $\alpha = 1$, we use respectively $v_8'(x) = 2x(1 + x^2)^{-1}$ and $v_8''(x) = (2 - 2x^2)(1 + x^2)^{-2}$, instead of \eqref{e:dav8}. Table \ref{t:testv8} shows the errors for six different approximations of $\partial_x\mathcal D^\alpha$, taking $N = 256$, $L = 30.4$. Even though the accuracy is rather low, it is nevertheless enough to get a rough idea of the shape of $\partial_x\mathcal D^\alpha v_8$.
\begin{table}[htb!]
\centering
\begin{tabular}{|c|c|}
\cline{2-2} \multicolumn{1}{c|}{} & $v_8(x)$
    \cr
\hline $\mathbf M^{(1)}$ &$2.1199\cdot10^{-2}$
    \cr
\hline $\mathbf M^{(1,2)}$ &$5.7432\cdot10^{-3}$
    \cr
\hline $\mathbf M^{(1,2,3)}$ &$9.6408\cdot10^{-4}$
    \cr
\hline $\mathbf M^{(1,2,3,4)}$ &$9.3004\cdot10^{-4}$
    \cr
\hline $\mathbf M^{(1,2,3,4,5)}$ &$9.3002\cdot10^{-4}$
    \cr
\hline $\mathbf M^{(1,2,3,4,5,6)}$ &$9.3002\cdot10^{-4}$
    \cr
\hline
\end{tabular}
\caption{Given $v_8(x) = \log(1 + x^2)$, we have computed its fractional derivative for six different-order approximations of $\partial_x\mathcal D^\alpha$, for $\alpha = j / 1000$, $j = 0, \ldots, 1000$. The column shows $\max_{\alpha\in[0,1]}E_8^{(1)}(\alpha)$, $\max_{\alpha\in[0,1]}E_8^{(1,2)}(\alpha)$, etc. The experiments have been done with $L = 30.4$ and $N = 256$.} \label{t:testv8}
\end{table}

Summarizing, although the method is best suited for regular functions, it gives acceptable results for functions far less regular, and it can even be of some utility when applied to unbounded functions with \emph{slow} growth.

\subsection{Discretization in time of the evolution problem \eqref{e:uoverR}}

\label{s:timediscretization}

The structure of \eqref{e:uoverR} suggests using a so-called implicit-explicit (IMEX) scheme, where the leading term is treated implicitly and the other terms are treated explicitly. In this paper, we have chosen the quite popular second-order semi-implicit backward differentiation formula (SBDF) \cite{ascher}, which, when applied to \eqref{e:uoverR}, takes the following form:
\begin{equation}
\label{e:SBDF2}
\begin{split}
\left(\frac{3}{2} - \Delta t\tau\partial_x^3\right) v^{(n+1)} =  2v^{(n)} & - \frac{1}{2}v^{(n-1)} + 2\Delta t\left[ \partial_x\mathcal D^\alpha v^{(n)} - \partial_x((v^{(n)})^2)\right]
    \cr
& - \Delta t\left[\partial_x\mathcal D^\alpha v^{(n-1)} - \partial_x((v^{(n-1)})^2)\right],
\end{split}
\end{equation}

\noindent where $v^{(n)}$ denotes the approximation of $v$ at time $t^{(n)} = n\Delta t$. We approximate the fractional derivative terms by means of $[\partial_x\mathcal D^\alpha]^{(1,2,3,4,5,6)}$; and the first-order derivative $\partial_x$ via \eqref{e:u_xu_s}. Once the right-hand side of \eqref{e:SBDF2} has been obtained, we observe that, from \eqref{e:udiffexp}, $(3 / 2 - \Delta\tau\partial_x^3)$ can be seen in the Fourier side as a very sparse matrix. Therefore, in view of the Fourier representation \eqref{e:approxu(s)}, it is possible to compute the $2N$ Fourier coefficients $\hat u^{(n+1)}(k)$, $-N\le k\le N-1$, in just $\mathcal O(N)$ operations. Observe also that we need two initial values to initialize the scheme \eqref{e:SBDF2}, $v^{(0)}$ and $v^{(1)}$. We take as $v^{(0)}$ an even extension of the initial data, while $v^{(1)}$ is obtained by a first-order semi-implicit Euler scheme,
\begin{equation*}
\left(1 - \Delta t\tau\partial_x^3\right) v^{(n+1)} =  v^{(n)}  + \Delta t\left[\partial_x\mathcal D^\alpha v^{(n)} -\partial_x((v^{(n)})^2)\right],
\end{equation*}

\noindent together with fourth-order Richardson extrapolation, in a similar way to, for instance, \eqref{e:extrapolation1}.

Let us make a couple of observations. The first one is that, although $\partial_x\mathcal D^\alpha$ is a linear operator, it is however a non-local one. Therefore, even if it is of course possible to write $\partial_x\mathcal D^\alpha v^{(n+1)}$ implicitly, it requires inverting a very dense matrix at each time-step, which is very inefficient. The second observation is that we have also tested the higher-order IMEX schemes in \cite{ascher}, but they appear to require bigger stability constraints on $\Delta t$. On the other hand, \eqref{e:SBDF2} is very stable. In our numerical experiments, unless otherwise indicated, we have taken $\Delta t = 0.01$, which is a relatively large value, but that allows to reach very long times with reasonable accuracy; smaller values of $\Delta t$ give no significant improvement in fact.

\section{Numerical experiments}

\label{s:experiments}
As an illustration of the method presented here, we give examples in which the initial condition induces the appearance of solutions that approach traveling wave solutions of \eqref{e:uoverR} as $t$ becomes large.

We recall that traveling wave solutions are solutions of the form $v(x,t)=\phi(\xi)$, where $\xi=x-ct$, for some constant wave speed $c$; and $\phi(\xi)$ approaches constant values as $\xi$ tends to $+\infty$ and to $-\infty$. Namely, if $\phi$ depends on $x$ and $t$ only through the traveling wave variable $\xi$, then ${\cal D}^{\alpha}\phi$ does too, and $\phi$ must satisfy, integrating once with respect to $\xi$,
\begin{equation}\label{TWP}
h(\phi)= {\cal D}^\alpha \phi+\tau\phi'',
\quad h(\phi) := -c(\phi-\phi_-)+\phi^2-\phi_-^2 ,
\end{equation}

\noindent where
\begin{equation*}
   c=\phi_+ + \phi_-,
\end{equation*}

\noindent and $\phi_-$ and $\phi_+$ are constants such that
\begin{equation*}
\lim_{\xi\to-\infty}\phi(\xi) = \phi_-,\qquad \lim_{\xi\to +\infty}\phi(\xi) = \phi_+,
\end{equation*}

\noindent and such that the entropy condition
\begin{equation*}
  \phi_- > \phi_+
\end{equation*}

\noindent is satisfied. For definiteness, we take
\begin{equation*}
\phi_-=1 ,\quad \phi_+=0 \quad \mbox{(i.e. $c=1$)},
\end{equation*}

\noindent and, hence, an initial condition that satisfies these far-field values is needed.

Traveling wave solutions have been studied in \cite{AHS} with $\tau=0$, and in \cite{ACH} with $\tau>0$. We recall that traveling waves for $\tau=0$ are monotone, as it is the case for the classical (or local) Burgers equation (see \cite{BonaSchon}). In \cite{ACH}, the existence of traveling waves for all values of $\tau$ is shown, as well as the dynamic stability of traveling waves, provided that they are monotone. The analysis guarantees the monotonicity of the waves for small values of the parameter $\tau$. In addition, we were able to show that the monotonicity in the {\it tail} of the waves is dominated by an algebraic decay of the form $\xi^{-\alpha}$, as $\xi \to +\infty$, by using results of Fractional Calculus that can be applied to the equation linearized in the tail (see \cite{GM2}). We do not have a quantification of the critical $\tau$ for which the transition from monotone to oscillatory behavior occurs, however. The mechanism of oscillations depends on $\alpha$ and on $\tau$; values of $\alpha$ close to $1$ damp the oscillations generated by the third order term at a larger values of $\tau$ than do smaller values of $\alpha$. And, as for the local case, the larger the value of $\tau$, the larger the amplitude of the oscillations for any fixed $\alpha\in(0,1)$.

The numerical examples performed in \cite{cuesta2014} by means of finite difference methods aim to complement the rigorous study and focus on the traveling wave equation (\ref{TWP}). The method developed in the current paper enables us to compute solutions for very long times with a relatively small number of nodes, and thus illustrate the convergence of solutions to traveling waves even for very large values of $\tau$. To illustrate this, we have simulated \eqref{e:uoverR}, supplemented with the initial datum
\begin{equation}
\label{e:u0}
v(x, 0) = \frac{1 - \tanh(x)}{2},
\end{equation}

\noindent i.e., a function with exponential decay. After applying the change of variable \eqref{e:changeofvariable} to the initial datum and getting a function defined over $[0,\pi]$, we have considered again an even extension, which is enough for our purposes, in order to get a function defined over the whole period. Let us mention, however, that there are extensions of a function defined over half a period that are smoother than the even extension (see v.g. \cite{morton2009} and \cite{hybrechs2010}), i.e., they make a more efficient use of the Fourier modes, which decay more quickly. Anyway, this problem is only posed in the initial data, since we do our simulations always over the full period.

In order to test \eqref{e:SBDF2}, we have taken $\alpha = 1/3$, $\tau = 1$, $N = 128$, and $L = 20$. We have done simulations with different time steps $\Delta t$, until $t = 20$; thus, for a given $\Delta t$, $v_{num}(x,20,\Delta t)$ denotes the numerical approximation of $v(x,20)$ using $\Delta t$. In order to apply a formula like \eqref{e:mualpha}, we would need the explicit exact solution of $v(x, 20)$, which is not available; instead, we use the $v_{num}(x,20,\Delta t)$ obtained for a very small value of $\Delta t$, say, $\Delta t = 10^{-4}$. Thus, we approximate the error in discrete $L^\infty$-norm as
\begin{equation*}
E(\Delta t) \equiv \|v_{num}(x, 20, \Delta t) - v_{num}(x, 20, \Delta t = 10^{-4})\|_{\infty}.
\end{equation*}

\noindent Table \ref{e:errorsevol} shows the evolution of $E(\Delta t)$ for different values of $\Delta t$. On the one hand, the errors are very small; on the other hand, the columns with $\log_2(E(\Delta t) / E(\Delta t/2))$ confirms clearly the second order of \eqref{e:SBDF2}.
\begin{table}[htb!]
    \centering
\begin{tabular}{|l|c|c|}
\hline \multicolumn{1}{|c|}{\phantom{\bigg(}$\Delta t$\phantom{\bigg)}} & $E(\Delta t)$ & $\log_2\left(\frac{E(\Delta t)}{E(\Delta t / 2)}\right)$
    \cr
\hline $1/5$ & $1.4801\cdot10^{-3}$ & $1.9704$
    \cr
\hline $1/10$ & $3.7769\cdot10^{-4}$ & $1.9941$
    \cr
\hline $1/20$ & $9.4811\cdot10^{-5}$ & $1.9985$
    \cr
\hline $1/40$ & $2.3727\cdot10^{-5}$ & $1.9997$
    \cr
\hline $1/80$ & $5.9332\cdot10^{-6}$ & $2.0002$
    \cr
\hline $1/160$ & $1.4831\cdot10^{-6}$ & $2.0012$
    \cr
\hline $1/320$ & $3.7048\cdot10^{-7}$
    \cr
\cline{1-2}
\end{tabular}
\caption{Errors in discrete $L^\infty$-norm of $v_{num}$ at $t = 20$, for different $\Delta t$. Since the exact solution is not available, we use instead $v_{num}(x, 20, \Delta t = 10^{-4})$ as reference. The results confirm the quadratic convergence rate of \eqref{e:SBDF2}.}
\label{e:errorsevol}
\end{table}

\begin{figure}[htb!]
\centering
\subfigure[Solutions against $x$.]
{
\includegraphics[height=0.37\textwidth,clip=true]{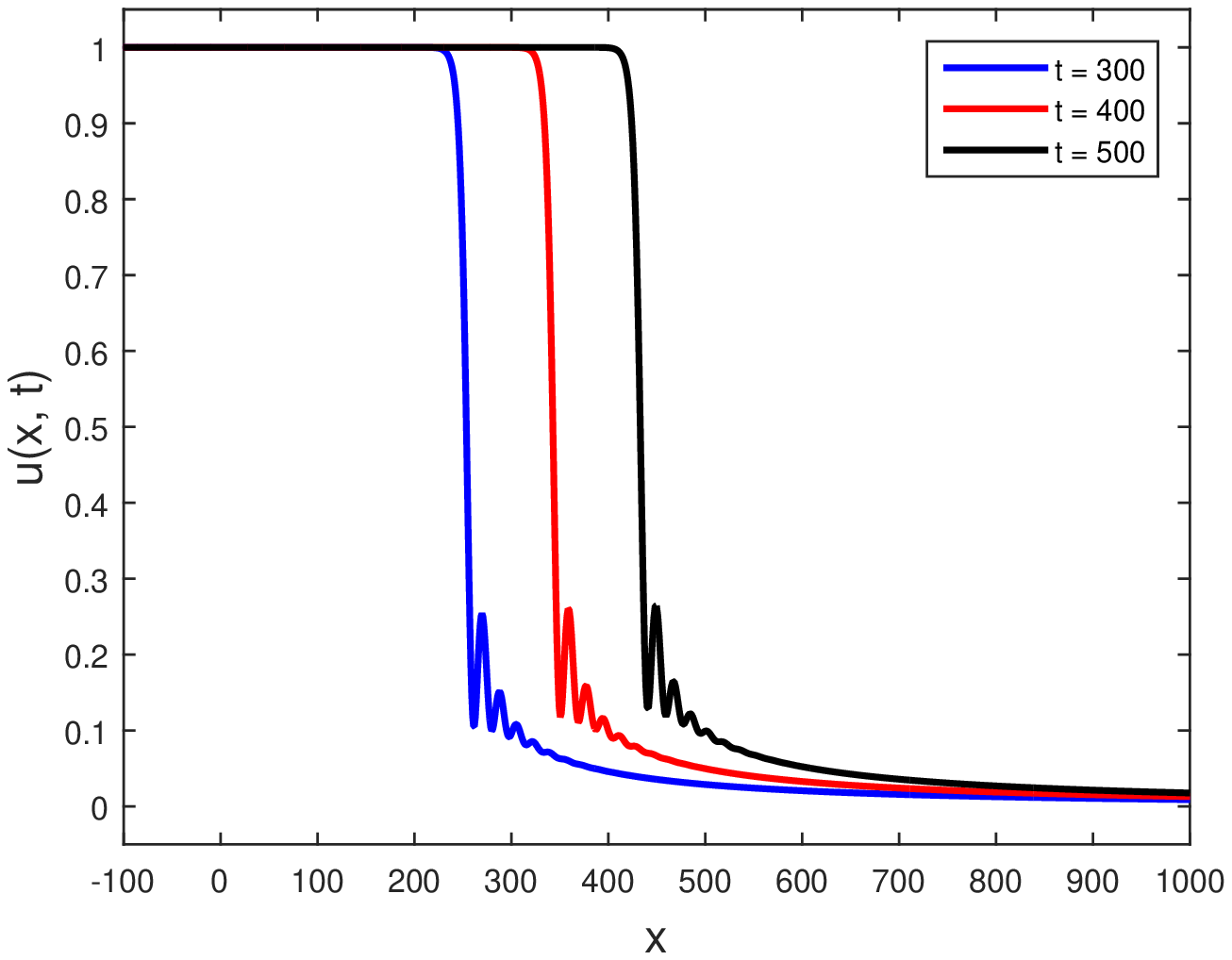}
\label{TWstau10a13}
}
\subfigure[Solutions against $x-ct$ with $c=1$.]
{
\includegraphics[height=0.37\textwidth,clip=true]{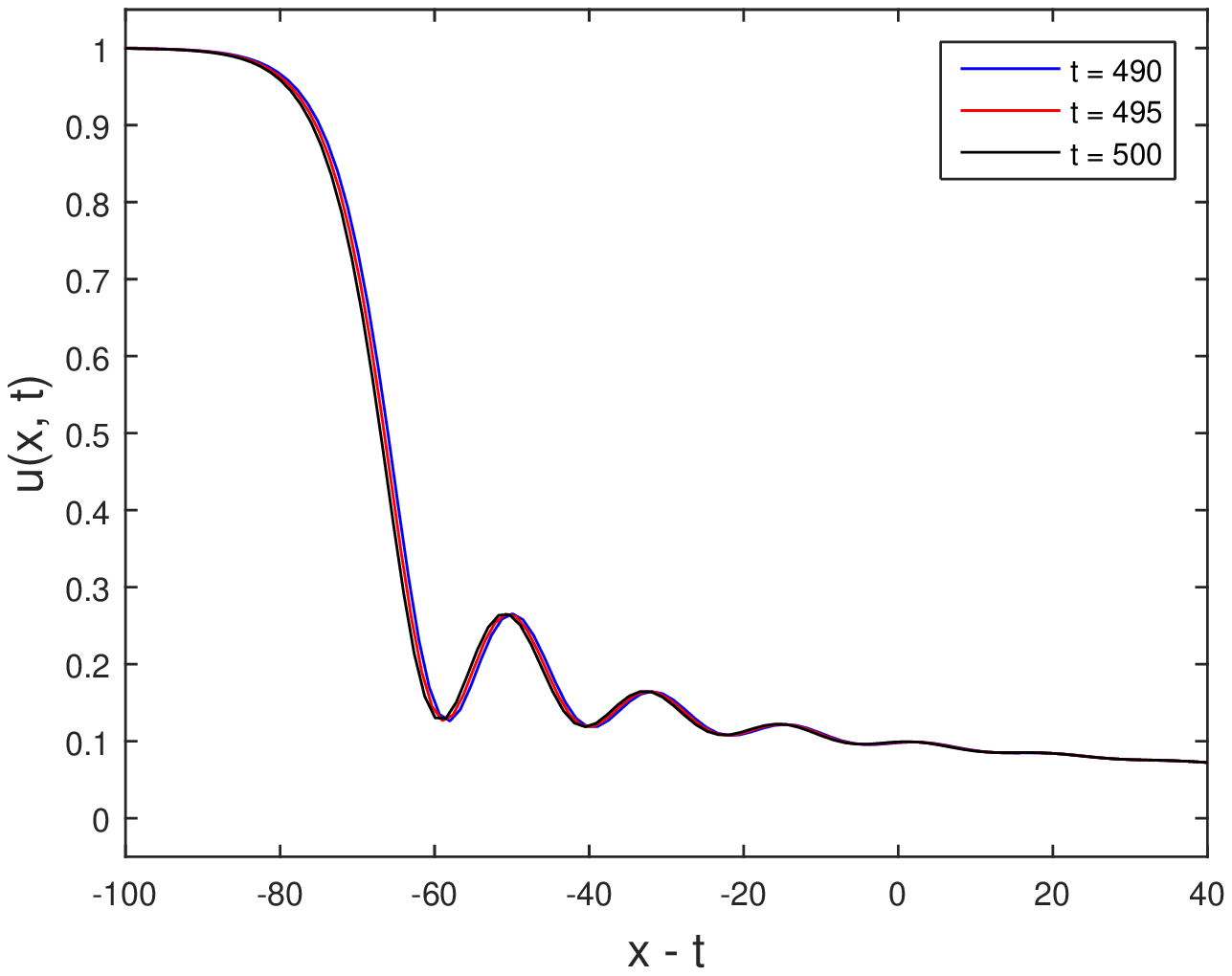}
\label{TWstau10a13II}
}
\caption{
Results for $\alpha=1/3$ with $\tau=10$. Here $L=500$, $\Delta t=0.01$ and $N=2048$.  \subref{TWstau10a13} shows solution profiles at $t=100$ (solid blue), $t=300$ (solid red), and $t=500$ (solid black).
\subref{TWstau10a13II} shows solution profiles against $x-ct$, with $c=1$, at $t=490$ (blue), $t=495$ (red), and $t=500$ (black).
}
\label{evoltau10a23}
\end{figure}

In Figure \ref{evoltau10a23}, we have taken $\alpha = 1/3$, $\tau = 10$, $N = 2048$. As mentioned earlier, the optimal choice of $L$ can be a delicate question. In our case, in order to simulate for very long times, a sufficiently large value of $L$ is required, if we want to have a decent resolution. In this example, $L = 500$ is enough to reach as far as $t = 500$, with $\Delta t = 10^{-2}$. Let us recall, however, that we are using a spectral method and that spatial resolution is not as important as having a large enough number of Fourier modes. Indeed, if the Fourier decomposition \eqref{e:approxu(s)} is known, and the values of $\hat u(k)$ decay fast enough as $|k|$ grows, it is straightforward to evaluate $u$ at any additional points, getting spectrally accurate values. Furthermore, known \eqref{e:approxu(s)}, it is even possible to modify $L$ and the number of nodes at any given time step, if needed.

\begin{figure}[htb!]
\centering
\includegraphics[width=0.48\textwidth,clip=true]{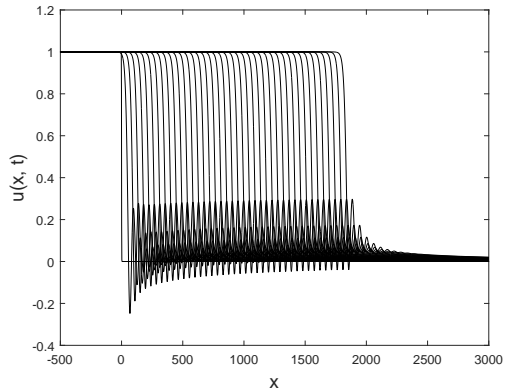}
\caption{
Results for $\alpha = 1 / 3$ with $\tau = 100$. Here $L = 2000$, $\Delta t = 0.01$ and $N = 4096$. We have plotted the results at times $t = 0, 50, \ldots, 2000$.
}
\label{f:evol4096}
\end{figure}

Finally, in Figure \ref{f:evol4096}, we have taken $\alpha = 1/3$, $\tau = 100$, $N = 4096$, with \eqref{e:u0} again as initial datum. For bigger $\tau$, it is necessary to consider longer times in order to see the stabilization to a traveling wave solution; hence, we have computed the evolution until $t = 2000$. We finally remark that the highest Fourier modes are of the order of $\mathcal O(10^{-13})$, which means that we could continue the simulation for a much longer time.

%\section{Conclusions}

%\label{s:conclusions}

%In this paper, we have developed a new pseudo-spectral method to solve the initial value problem associated to the non-local KdV-Burgers equation \eqref{e:uoverR}. This method allows us to simulate accurately \eqref{e:uoverR} until extremely large times, e.g., $t = 2000$, and is thus especially convenient when the solutions approach traveling wave solutions of \eqref{e:uoverR}, as $t$ becomes large.

%The central part of this paper is devoted to the accurate numerical computation of $\partial_x\mathcal D^\alpha$ on the whole real line, for which \eqref{e:assympt} plays a fundamental r\^ole. Even if \eqref{e:assympt} has been obtained empirically, its adequacy is confirmed by a number of tests involving functions with different decays and regularity. Nevertheless, it would be interesting to offer an analytical proof for it, although this lies beyond the scope of this paper.

\section*{Acknowledgments}

The authors want to thank the anonymous referees for their useful comments that have lead to important improvements of this paper.

The authors were supported by The Basque Government through the project IT641-13. They also acknowledge the support from the Spanish Ministry of Economy and Competitiveness through the projects MTM2011-24054 (F. de la Hoz), MTM2011-24109 (C. M. Cuesta) and MTM2014-53145-P (F. de la Hoz and C. M. Cuesta).

%\section*{References}

%\bibliographystyle{elsarticle-num}
%\bibliography{carlotapatxi}

\bibliographystyle{abbrv}

\def\cprime{$'$}

\end{document}